\documentclass[12pt]{article}
\usepackage{theorem}
\usepackage{amssymb}
\usepackage{graphicx}
\usepackage[mathscr]{eucal}
\usepackage{amsbsy}
\usepackage{amsmath}
\newtheorem{prop}{}[section]

{\theorembodyfont{\upshape} \newtheorem{rema}[prop]{}}
\newcommand{\boma}[1]{{\mbox{\boldmath $#1$} }}

\begin{document}
\newcommand{\uper}[1]{\stackrel{\barray{c} {~} \\ \mbox{\footnotesize{#1}}\farray}{\longrightarrow} }
\newcommand{\nop}[1]{ \|#1\|_{\piu} }
\newcommand{\no}[1]{ \|#1\| }
\newcommand{\nom}[1]{ \|#1\|_{\meno} }
\newcommand{\uu}[1]{e^{#1 \AA}}
\newcommand{\UD}[1]{e^{#1 \Delta}}
\newcommand{\bb}[1]{\mathbb{{#1}}}
\newcommand{\HO}[1]{\bb{H}^{{#1}}}
\newcommand{\Hz}[1]{\bb{H}^{{#1}}_{\zz}}
\newcommand{\Hs}[1]{\bb{H}^{{#1}}_{\ss}}
\newcommand{\Hg}[1]{\bb{H}^{{#1}}_{\gg}}
\newcommand{\HM}[1]{\bb{H}^{{#1}}_{\so}}
\newcommand{\hz}[1]{H^{{#1}}_{\zz}}
\def\ug{u^G}
\def\gaamma{\widehat{\gamma}}
\def\Rc{R_{\star}}
\def\Ddd{\widehat{\Dd}}
\def\Rrr{\widehat{\Rr}}
\def\Tcc{\widehat{T}_{\tt{c}}}
\def\eep{\widehat{\epsilon}}
\def\mum{\hat{\mu}}
\def\uz{u_{*}}
\def\R{{R}\,}
\def\ti{{\tt{t}}}
\def\ef{\psi}
\def\fun{\mathcal{F}}
\def\fun{{\tt f}}
\def\tvainf{\vspace{-0.4cm} \barray{ccc} \vspace{-0,1cm}{~}
\\ \vspace{-0.2cm} \longrightarrow \\ \vspace{-0.2cm} \scriptstyle{T \vain + \infty} \farray}
\def\De{F}
\def\er{\epsilon}
\def\erd{\er_0}
\def\Tn{T_{\star}}
\def\Tc{T_{\tt{c}}}
\def\Tb{T_{\tt{b}}}
\def\Tl{\mathscr{T}}
\def\Tm{T}
\def\Ta{T_{\tt{a}}}
\def\ua{u_{\tt{a}}}
\def\Tg{T_{G}}
\def\Tgg{T_{I}}
\def\Tw{T_{w}}
\def\Ts{T_{\Ss}}
\def\Tr{\Tl}
\def\Sp{\Ss'}
\def\Tsp{T_{\Sp}}
\def\vsm{\vspace{-0.1cm}\noindent}
\def\comple{\scriptscriptstyle{\complessi}}
\def\nume{0.407}
\def\numerob{0.00724}
\def\deln{7/10}
\def\delnn{\dd{7 \over 10}}
\def\e{c}
\def\p{p}
\def\z{z}
\def\symd{{\mathfrak S}_d}
\def\del{\omega}
\def\Del{\delta}
\def\Di{\Delta}
\def\Ss{{\mathscr{S}}}
\def\Ww{{\mathscr{W}}}
\def\mmu{\hat{\mu}}
\def\rot{\mbox{rot}\,}
\def\curl{\mbox{curl}\,}
\def\Mm{\mathscr M}
\def\XS{\boma{x}}
\def\TS{\boma{t}}
\def\Lam{\boma{\eta}}
\def\DS{\boma{\rho}}
\def\KS{\boma{k}}
\def\LS{\boma{\lambda}}
\def\PR{\boma{p}}
\def\VS{\boma{v}}
\def\ski{\! \! \! \! \! \! \! \! \! \! \! \! \! \!}
\def\h{L}
\def\EM{M}
\def\EMP{M'}
\def\Rr{{\mathscr{R}}}
\def\Zz{{\mathscr{Z}}}
\def\E{E}
\def\FFf{\mathscr{F}}
\def\A{F}
\def\Xim{\Xi_{\meno}}
\def\Ximn{\Xi_{n-1}}
\def\lan{\lambda}
\def\om{\omega}
\def\Om{\Omega}
\def\Sim{\Sigm}
\def\Sip{\Delta \Sigm}
\def\Sigm{{\mathscr{S}}}
\def\Ki{{\mathscr{K}}}
\def\Hi{{\mathscr{H}}}
\def\zz{{\scriptscriptstyle{0}}}
\def\ss{{\scriptscriptstyle{\Sigma}}}
\def\gg{{\scriptscriptstyle{\Gamma}}}
\def\so{\ss \zz}
\def\Dv{\bb{\DD}'}
\def\Dz{\bb{\DD}'_{\zz}}
\def\Ds{\bb{\DD}'_{\ss}}
\def\Dsz{\bb{\DD}'_{\so}}
\def\Dg{\bb{\DD}'_{\gg}}
\def\Ls{\bb{L}^2_{\ss}}
\def\Lg{\bb{L}^2_{\gg}}
\def\bF{{\bb{V}}}
\def\Fz{\bF_{\zz}}
\def\Fs{\bF_\ss}
\def\Fg{\bF_\gg}
\def\Pre{P}
\def\UU{{\mathcal U}}
\def\fiapp{\phi}
\def\PU{P1}
\def\PD{P2}
\def\PT{P3}
\def\PQ{P4}
\def\PC{P5}
\def\PS{P6}
\def\Q{P6}
\def\X{Q2}
\def\Xp{Q3}
\def\Vi{V}
\def\bVi{\bb{V}}
\def\K{V}
\def\Ks{\bb{\K}_\ss}
\def\Kz{\bb{\K}_0}
\def\KM{\bb{\K}_{\, \so}}
\def\HGG{\bb{H}^\G}
\def\HG{\bb{H}^\G_{\so}}
\def\EG{{\mathfrak{P}}^{\G}}
\def\G{G}
\def\de{\delta}
\def\esp{\sigma}
\def\dd{\displaystyle}
\def\LP{\mathfrak{L}}
\def\dive{\mbox{div}}
\def\la{\langle}
\def\ra{\rangle}
\def\um{u_{\meno}}
\def\uv{\mu_{\meno}}
\def\Fp{ {\textbf F_{\piu}} }
\def\Ff{ {\textbf F} }
\def\Fm{ {\textbf F_{\meno}} }
\def\Eb{ {\textbf E} }
\def\piu{\scriptscriptstyle{+}}
\def\meno{\scriptscriptstyle{-}}
\def\omeno{\scriptscriptstyle{\ominus}}
\def\Tt{ {\mathscr T} }
\def\Xx{ {\textbf X} }
\def\Yy{ {\textbf Y} }
\def\Ee{ {\textbf E} }
\def\VP{{\mbox{\tt VP}}}
\def\CP{{\mbox{\tt CP}}}
\def\cp{$\CP(f_0, t_0)\,$}
\def\cop{$\CP(f_0)\,$}
\def\copn{$\CP_n(f_0)\,$}
\def\vp{$\VP(f_0, t_0)\,$}
\def\vop{$\VP(f_0)\,$}
\def\vopn{$\VP_n(f_0)\,$}
\def\vopdue{$\VP_2(f_0)\,$}
\def\leqs{\leqslant}
\def\geqs{\geqslant}
\def\mat{{\frak g}}
\def\tG{t_{\scriptscriptstyle{G}}}
\def\tN{t_{\scriptscriptstyle{N}}}
\def\TK{t_{\scriptscriptstyle{K}}}
\def\CK{C_{\scriptscriptstyle{K}}}
\def\CN{C_{\scriptscriptstyle{N}}}
\def\CG{C_{\scriptscriptstyle{G}}}
\def\CCG{{\mathscr{C}}_{\scriptscriptstyle{G}}}
\def\tf{{\tt f}}
\def\ta{{\tt a}}
\def\tc{{\tt c}}
\def\tF{{\tt R}}
\def\C{{\mathscr C}}
\def\P{{\mathscr P}}
\def\V{{\mathscr V}}
\def\TI{\tilde{I}}
\def\TJ{\tilde{J}}
\def\Lin{\mbox{Lin}}
\def\Hinfc{ H^{\infty}(\reali^d, \complessi) }
\def\Hnc{ H^{n}(\reali^d, \complessi) }
\def\Hmc{ H^{m}(\reali^d, \complessi) }
\def\Hac{ H^{a}(\reali^d, \complessi) }
\def\Dc{\DD(\reali^d, \complessi)}
\def\Dpc{\DD'(\reali^d, \complessi)}
\def\Sc{\SS(\reali^d, \complessi)}
\def\Spc{\SS'(\reali^d, \complessi)}
\def\Ldc{L^{2}(\reali^d, \complessi)}
\def\Lpc{L^{p}(\reali^d, \complessi)}
\def\Lqc{L^{q}(\reali^d, \complessi)}
\def\Lrc{L^{r}(\reali^d, \complessi)}
\def\Hinfr{ H^{\infty}(\reali^d, \reali) }
\def\Hnr{ H^{n}(\reali^d, \reali) }
\def\Hmr{ H^{m}(\reali^d, \reali) }
\def\Har{ H^{a}(\reali^d, \reali) }
\def\Dr{\DD(\reali^d, \reali)}
\def\Dpr{\DD'(\reali^d, \reali)}
\def\Sr{\SS(\reali^d, \reali)}
\def\Spr{\SS'(\reali^d, \reali)}
\def\Ldr{L^{2}(\reali^d, \reali)}
\def\Hinfk{ H^{\infty}(\reali^d, \KKK) }
\def\Hnk{ H^{n}(\reali^d, \KKK) }
\def\Hmk{ H^{m}(\reali^d, \KKK) }
\def\Hak{ H^{a}(\reali^d, \KKK) }
\def\Dk{\DD(\reali^d, \KKK)}
\def\Dpk{\DD'(\reali^d, \KKK)}
\def\Sk{\SS(\reali^d, \KKK)}
\def\Spk{\SS'(\reali^d, \KKK)}
\def\Ldk{L^{2}(\reali^d, \KKK)}
\def\Knb{K^{best}_n}
\def\sc{\cdot}
\def\k{\mbox{{\tt k}}}
\def\x{\mbox{{\tt x}}}
\def\g{ {\textbf g} }
\def\QQQ{ {\textbf Q} }
\def\AAA{ {\textbf A} }
\def\gr{\mbox{gr}}
\def\sgr{\mbox{sgr}}
\def\loc{\mbox{loc}}
\def\PZ{{\Lambda}}
\def\PZAL{\mbox{P}^{0}_\alpha}
\def\epsilona{\epsilon^{\scriptscriptstyle{<}}}
\def\epsilonb{\epsilon^{\scriptscriptstyle{>}}}
\def\lgraffa{ \mbox{\Large $\{$ } \hskip -0.2cm}
\def\rgraffa{ \mbox{\Large $\}$ } }
\def\restriction{\upharpoonright}
\def\M{{\scriptscriptstyle{M}}}
\def\m{m}
\def\Fre{Fr\'echet~}
\def\I{{\mathcal N}}
\def\ap{{\scriptscriptstyle{ap}}}
\def\fiap{\varphi_{\ap}}
\def\dfiap{{\dot \varphi}_{\ap}}
\def\DDD{ {\mathfrak D} }
\def\BBB{ {\textbf B} }
\def\EEE{ {\textbf E} }
\def\GGG{ {\textbf G} }
\def\TTT{ {\textbf T} }
\def\KKK{ {\textbf K} }
\def\HHH{ {\textbf K} }
\def\FFi{ {\bf \Phi} }
\def\GGam{ {\bf \Gamma} }
\def\sc{ {\scriptstyle{\bullet} }}
\def\a{a}
\def\ep{\epsilon}
\def\c{\kappa}
\def\parn{\par \noindent}
\def\teta{M}
\def\elle{L}
\def\ro{\rho}
\def\al{\alpha}
\def\si{\sigma}
\def\be{\beta}
\def\ga{\gamma}
\def\te{\vartheta}
\def\ch{\chi}
\def\et{\eta}
\def\complessi{{\bf C}}
\def\len{{\bf L}}
\def\reali{{\bf R}}
\def\interi{{\bf Z}}
\def\Z{{\bf Z}}
\def\naturali{{\bf N}}
\def\Sfe{ {\bf S} }
\def\To{ {\bf T} }
\def\Td{ {\To}^d }
\def\Tt{ {\To}^3 }
\def\Zd{ \interi^d }
\def\Zt{ \interi^3 }
\def\Zet{{\mathscr{Z}}}
\def\Ze{\Zet^d}
\def\T1{{\textbf To}^{1}}
\def\es{s}
\def\ee{{E}}
\def\FF{\mathcal F}
\def\FFu{ {\textbf F_{1}} }
\def\FFd{ {\textbf F_{2}} }
\def\GG{{\mathcal G} }
\def\EE{{\mathcal E}}
\def\KK{{\mathcal K}}
\def\PP{{\mathcal P}}
\def\PPP{{\mathscr P}}
\def\PN{{\mathcal P}}
\def\PPN{{\mathscr P}}
\def\QQ{{\mathcal Q}}
\def\J{J}
\def\Np{{\hat{N}}}
\def\Lp{{\hat{L}}}
\def\Jp{{\hat{J}}}
\def\Pp{{\hat{P}}}
\def\Pip{{\hat{\Pi}}}
\def\Vp{{\hat{V}}}
\def\Ep{{\hat{E}}}
\def\Gp{{\hat{G}}}
\def\Kp{{\hat{K}}}
\def\Ip{{\hat{I}}}
\def\Tp{{\hat{T}}}
\def\Mp{{\hat{M}}}
\def\La{\Lambda}
\def\Ga{\Gamma}
\def\Si{\Sigma}
\def\Upsi{\Upsilon}
\def\Gam{\Gamma}
\def\Gag{{\check{\Gamma}}}
\def\Lap{{\hat{\Lambda}}}
\def\Upsig{{\check{\Upsilon}}}
\def\Kg{{\check{K}}}
\def\ellp{{\hat{\ell}}}
\def\j{j}
\def\jp{{\hat{j}}}
\def\BB{{\mathcal B}}
\def\LL{{\mathcal L}}
\def\MM{{\mathcal U}}
\def\SS{{\mathcal S}}
\def\DD{D}
\def\Dd{{\mathcal D}}
\def\VV{{\mathcal V}}
\def\WW{{\mathcal W}}
\def\OO{{\mathcal O}}
\def\RR{{\mathcal R}}
\def\TT{{\mathcal T}}
\def\AA{{\mathcal A}}
\def\CC{{\mathcal C}}
\def\JJ{{\mathcal J}}
\def\NN{{\mathcal N}}
\def\HH{{\mathcal H}}
\def\XX{{\mathcal X}}
\def\XXX{{\mathscr X}}
\def\YY{{\mathcal Y}}
\def\ZZ{{\mathcal Z}}
\def\CC{{\mathcal C}}
\def\cir{{\scriptscriptstyle \circ}}
\def\circa{\thickapprox}
\def\vain{\rightarrow}
\def\salto{\vskip 0.2truecm \noindent}
\def\spazio{\vskip 0.5truecm \noindent}
\def\vs1{\vskip 1cm \noindent}
\def\fine{\hfill $\square$ \vskip 0.2cm \noindent}
\def\ffine{\hfill $\lozenge$ \vskip 0.2cm \noindent}
\newcommand{\rref}[1]{(\ref{#1})}
\def\beq{\begin{equation}}
\def\feq{\end{equation}}
\def\beqq{\begin{eqnarray}}
\def\feqq{\end{eqnarray}}
\def\barray{\begin{array}}
\def\farray{\end{array}}
%%%%%%%%% THIS NUMBERS EQUATIONS BY SECTIONS %%%%%%%%%%%%%
\makeatletter \@addtoreset{equation}{section}
\renewcommand{\theequation}{\thesection.\arabic{equation}}
%\thesection instead of \arabic{section} for correct equation numbering
% in appendices
\makeatother
%%%%%%%%%%%%%%%%%%%%%%%%%%%INTESTAZIONE%%%%%%%%%%%%%%%%%%%%%%%%%%%%%%%
\begin{titlepage}
{~}
\vspace{-2cm}
\begin{center}
{\huge On the Reynolds number expansion for the Navier-Stokes equations}
\end{center}
\vspace{0.5truecm}
\begin{center}
{\large
Carlo Morosi$\,{}^a$, Livio Pizzocchero$\,{}^b$({\footnote{Corresponding author}})} \\
\vspace{0.5truecm} ${}^a$ Dipartimento di Matematica, Politecnico di Milano,
\\ P.za L. da Vinci 32, I-20133 Milano, Italy \\
e--mail: carlo.morosi@polimi.it \\
${}^b$ Dipartimento di Matematica, Universit\`a di Milano\\
Via C. Saldini 50, I-20133 Milano, Italy\\
and Istituto Nazionale di Fisica Nucleare, Sezione di Milano, Italy \\
e--mail: livio.pizzocchero@unimi.it
\end{center}
\begin{abstract}
In a previous paper of ours \cite{appeul} we have considered
the incompressible Navier-Stokes (NS) equations on a $d$-dimensional torus
$\Td$, in the functional setting of the Sobolev spaces $\HM{n}(\Td)$
of divergence free, zero mean vector fields ($n > d/2+1)$.
In the cited work we have presented a general setting for the \textsl{a posteriori}
analysis of approximate solutions of the NS
Cauchy problem; given any approximate solution $\ua$,
this allows to infer a lower bound $\Tc$ on the time of existence
of the exact solution $u$ and to construct a
function $\Rr_n$ such that $\| u(t) - \ua(t) \|_n \leqs \Rr_n(t)$ for all $t \in [0,\Tc)$.
In certain cases it is $\Tc = + \infty$, so global existence is granted
for $u$. In the present paper the framework of \cite{appeul} is applied using as
an approximate solution an expansion $u^N(t) = \sum_{j=0}^N R^j u_j(t)$,
where $R$ is the Reynolds number. This allows, amongst else, to derive
the global existence of $u$ when $R$ is below some critical value
$R_{*}$ (increasing with $N$ in the examples
that we analyze). After a general discussion about the Reynolds expansion
and its \textsl{a posteriori} analysis, we consider the expansions of orders $N=1,2,5$
in dimension $d=3$, with the initial datum of Behr, Ne$\check{\mbox{c}}$as and Wu
\cite{Nec}. Computations of order $N=5$ yield a quantitative improvement of
the results previously obtained for this initial datum
in \cite{appeul}, where a Galerkin approximate solution was employed in place
of the Reynolds expansion.
\end{abstract}
\vspace{0.2cm} \noindent
\hfill \par
\par \vspace{0.05truecm} \noindent \textbf{AMS 2000 Subject classifications:} 35Q30, 76D03, 76D05.
\end{titlepage}
\section{Introduction}
\label{intro}
The incompressible Navier-Stokes (NS) equations with no external forces and periodic
boundary conditions can be written as
\beq {\partial u \over \partial \ti}  = \nu \Delta u + \PPP(u,u)~,
\label{eulnu} \feq
where: $\nu \in (0,+\infty)$ is the viscosity coefficient;
$u= u(x, \ti)$ is the divergence free velocity field;
the space variables $x = (x_s)_{s=1,...,d}$ belong to the torus
$\Td$
(and yield the derivatives $\partial_s := \partial/\partial x_s$);
$\Delta := \sum_{s=1}^d \partial_{s s}$ is the Laplacian. Furthermore,
$\PPP$ is the bilinear map defined as follows:
for all sufficiently regular velocity fields $v, w$ on $\Td$,
\beq \PPP(v, w) := - \LP(v \sc \partial w) \label{ppp} \feq
where $(v \sc \partial w)_r := \sum_{s=1}^d v_s \partial_s w_r$ ($r=1,...,d$),
and $\LP$ is the Leray projection onto the space of divergence free vector fields.
The dimension $d$ is arbitrary in the general setting of this paper,
but we put $d=3$ in the application of the last section.
%%From the physical viewpoint,
%%the velocity fields (at any fixed time) should be maps from $\Td$ to $\reali^d$; however,
%%we can harmlessly consider maps $\Td \vain \complessi^d$.
\par
Let us introduce the rescaled time $t$ and the Reynolds number $R$, setting
\beq t := \nu \ti~, \qquad R := {1 \over \nu}~; \label{setting} \feq
then Eq \rref{eulnu} takes the form
\beq {\partial u \over \partial t}  = \Delta u + R \, \PPP(u,u)~, \label{eul} \feq
that will be the standard of this paper.
Our functional setting for Eq. \rref{eul} is based on the Sobolev spaces
\beq \HM{n}(\Td) \equiv \HM{n} := \{ v : \Td \vain \reali^d~|
 \langle v \rangle = 0, \dive \, v = 0, \sqrt{-\Delta}^{\,n} v \in \mathbb{L}^2(\Td) \}  \feq
 (where $\langle~\rangle$ indicates the mean over $\Td$);
 for any real $n$, the above space is equipped with
the inner product $\la v | w \ra_n := \la \sqrt{-\Delta}^{\,n} v |  \sqrt{-\Delta}^{\,n} w \ra_{L^2}$
and with the corresponding norm $\| ~ \|_n$.
\par
In our paper \cite{appeul} we have outlined
a general framework to obtain quantitative
information on the exact solution $u$ of the NS
Cauchy problem analyzing \textsl{a posteriori} an
approximate solution. To be more precise,
consider the NS equation \rref{eul} with a specified
initial condition $u(x,0) = \uz(x)$; let $\ua : \Td \times [0,\Ta) \vain \reali^d$ be
an approximate solution of this Cauchy problem, and consider (for $n > d/2 + 1$) the
Sobolev norms
\parn
\vbox{
\beq \| \big({\partial \ua \over \partial t} - \Delta \ua
- R \, \PPP(\ua, \ua))(t)\|_n~,
\qquad \| \ua(0) - \uz \|_n~, \label{norms1} \feq
\beq \| \ua(t) \|_n~, \quad \| \ua(t) \|_{n+1},
\label{norms2} \feq
}
where $t \in [0, \Ta)$ and $\ua(t) := \ua(\cdot, t)$.
The norms in \rref{norms1} control the \emph{differential and datum
errors} of $\ua$, while the norms in \rref{norms2} refer to the growth
of $\ua$.
The approach of \cite{appeul} relies
on the so-called \emph{control inequalities};
these consist of a differential inequality and of
an inequality on the initial value,
determined by the norms \rref{norms1}
\rref{norms2} and involving an unknown function $\Rr_n : [0,\Tc) \vain [0,+\infty)$.
Assume the control inequalities to have
a solution $\Rr_n$, with a suitable domain $[0, \Tc)$;
then, according to \cite{appeul}, the solution $u$ of the NS equation \rref{eul} with initial datum $\uz$
exists (in a classical sense) on the time interval $[0,\Tc)$, and its distance from the
approximate solution admits the bound
\beq \| u(t) - \ua(t) \|_n \leqs \Rr_n(t) \qquad \mbox{for $t \in [0,\Tc)$}~. \label{buap} \feq
(For similar or related statements on the NS equations and other nonlinear
evolutionary PDEs, see the papers by Chernyshenko \textsl{et al.} \cite{Che},
Robinson \textsl{et al.} \cite{Rob} and our works \cite{uno} \cite{due} \cite{accau}). \par
In the present paper we apply the above framework choosing
as an approximate solution a polynomial in $R$ of
the form
\beq u^N(t) := \sum_{j=0}^{N} R^j u_j(t)~, \label{tay} \feq
where the terms $u_j(t)$ are determined requiring the
differential error to be $O(R^{N+1})$ for $R \vain 0$.
We emphasize that,
in our approach, the order $N$ could be large but is fixed;
so we are not considering the $N \vain + \infty$ limit,
i.e., the solution of the Cauchy problem for Eq. \rref{eul}
via a power series $u(t) = \sum_{j=0}^{+\infty} R^j u_j(t)$.
A theoretical analysis of the convergence issue
for such a series, in suitable function spaces, has
been developed by some authors, especially Cannone \cite{Can}
and Sinai \cite{Sin}. However, the approaches of these authors
yield convergence conditions (local or
global in time) which depend on the initial datum
$\uz$ only through its norm; moreover, these authors have dedicated little attention
to the strictly quantitative aspects of their analysis
(such as the evaluations of the constants in certain
inequalities). On the contrary,
our approach based on a finite order
approximant $u^N$ as in \rref{tay} has the following features. \par
\begin{itemize}
\item[(i)]
We produce estimates on the interval of existence
of the exact solution $u$ (and on its distance
from $u^N$) which depend on the fine structure of
the initial datum $\uz$ and not only on its norm;
the specific features of the initial datum
yielding these estimates are encoded in the
expression of the differential error $\partial u^N/\partial t - \Delta u^N
- R \, \PPP(u^N, u^N)$.
\item[(ii)]
Our analysis is fully quantitative: it relies
on explicit expressions for $u^N$ and its
errors and uses, amongst else, the estimates of \cite{cog} \cite{cok}
on the constants in certain inequalities about $\PPP$.
\end{itemize}
The above setting invites a computer assisted
approach: this can be readily set up when
the NS initial datum $\uz$ is sufficiently simple, say,
a Fourier polynomial. In this case
the approximant $u^N$ and the errors
$\partial u^N/\partial t - \Delta u^N - R \, \PPP(u^N, u^N)$,
$u^N(0) - \uz$, with their
Sobolev norms, can be determined via any package for
symbolic computation; after this a solution
$\Rr_n$ for the control inequalities can be
obtained numerically. More precisely, one
can try to satisfy them as equalities:
this amounts to solve the Cauchy problem
for a simple ODE in the unknown
function $t \mapsto \Rr_n(t)$. This ``control Cauchy problem''
is easily treated numerically. \par
In the present paper the above procedure is described in
general terms and then applied
with $d=3$ and $n=3$, choosing for
$\uz$ the so-called Behr-Ne$\check{\mbox{c}}$as-Wu
(BNW) initial datum \cite{Nec}. For the practical
implementation we use MATHEMATICA on a PC; first we work with
$N=1$ to introduce the method, and then
pass to the orders $N=2$, $N=5$. \par
In all the above cases,
the control Cauchy problem has a
solution $\Rr_3$ of domain $[0,+\infty)$
if $R$ is below some critical value $R_{*}$ (depending
on $N$); in this situation we can grant
global existence for the solution $u$
of the NS Cauchy problem, and
the inequality $\| u(t) - u^N(t) \|_3 \leqs \Rr_3(t)$
holds for all $t \in [0,+\infty)$.
For $R$ above $R_{*}$, the control
problem has a solution $\Rr_3$ on
a bounded interval $[0, \Tc)$; so,
the existence of the solution
$u$ of the NS Cauchy problem
and the inequality $\| u(t) - u^N(t) \|_3
\leqs \Rr_3(t)$ are granted only on this
interval. Passing from $N=1$ to $N=2$, and
from $N=2$ to $N=5$,
the critical value $R_{*}$ increases.
\par
In the final part of the paper the outcomes of the above Reynolds expansions are compared
with the results of \cite{appeul}, where
an approximate NS solution was constructed
for the BNW initial datum using the Galerkin method
with a set of 150 Fourier modes; it turns out that
this Galerkin approach is quantitatively
equivalent to the Reynolds expansion of
order $N=2$, while the expansion of order
$N=5$ gives much better results
concerning the global existence of $u$
and its distance from the approximate solution. \par
The paper is organized as follows. After fixing some
basic notations, Section \ref{prelim} reviews the general
setting of \cite{appeul} for approximate NS solutions
(in a reformulation suitable for our purposes,
where the Reynolds number $R =1/\nu$ and
the rescaled time $t = \nu \ti$ are
preferred to the variables $\nu, \ti$). Sections \ref{reynolds} and \ref{imple}
present the general Reynolds expansion
\rref{tay}, the control inequalities for it
and some basic computational rules for
the practical implementation of this approach.
Section \ref{behr} applies the previous
framework to the the BNW initial datum,
and makes a comparison with the Galerkin approach of
\cite{appeul} for the same datum. \par
In a forthcoming paper \cite{Forth} the Reynolds expansion and
the related control equation will be applied
to the BNW datum for larger values of $N$, and then employed for
other initial data of interest in the NS community,
namely, the ``vortices'' of Taylor-Green \cite{Tay} and Kida \cite{Kid} .
\section{Preliminaries}
\label{prelim}
Throughout the paper we fix a space
dimension $d \in \{2,3,...\}$; in the application of section \ref{behr} we will
put $d=3$.
For $a, b$ in $\reali^d$ or $\complessi^d$ we put $a \sc b :=
\sum_{r=1}^d a_r b_r$ and $|a| := \sqrt{\overline{a} \sc a}$,
where $\overline{\phantom{x}}$ denotes the complex conjugation. \par
Let us consider the $d$-dimensional torus $\Td :=
(\reali/2 \pi \interi)^d$;
the setting outlined hereafter for function spaces and
NS equations on the torus is similar to the one
of \cite{accau} \cite{appeul}.
%%A main difference between
%%these papers and the present one is that here we deal
%%with spaces of complex valued, rather than real valued
%%functions; some results in \cite{accau} \cite{appeul} (and
%%in other works of ours), originally formulated in a real
%%framework, hold as well in the present, complexified
%%setting. \par
In the sequel we employ the Fourier basis
made of the functions
({\footnote{In \cite{accau} \cite{appeul},
the normalization factor $(2 \pi)^{-d/2}$ was included in
the definition of $e_k$.}})
\beq e_k : \Td \vain \complessi~, \qquad e_k(x) :=
e^{i k \sc x} \qquad (k \in \Zd)~. \feq
Let us consider the space
$D'(\Td, \reali^d) \equiv D'$ of real distributions on $\Td$
%%i.e., the topological dual of $D := C^\infty(\Td, \complessi)$
({\footnote{This is defined in terms
of the space $D'(\Td, \complessi) \equiv D'_{\complessi}$
of complex distributions on $\Td$, which is the topological
dual of $C^{\infty}(\Td, \complessi)$. $D'_{\complessi}$ is
known to carry a complex
conjugation $\overline{\phantom{f}}$,
see, e.g., \cite{accau}; real
distributions $v$ on $\Td$ are characterized
by the condition $\overline{v}=v$.}}).
Any $v \in D'$ has a weakly convergent Fourier expansion
$v = \sum_{k \in \Zd} v_k e_k$, where
$v_k := (2 \pi)^{-d} \langle v,  e_{-k} \rangle $ (the notation indicates
the action of $v$ on the test function $e_{-k}$); of course the reality of
$v$ is expressed by the condition $\overline{v_k} = v_{-k}$.
\par
We will often be interested in the spaces
$L^p(\Td, \reali) \equiv L^p$ and, in particular,
in the real Hilbert space $L^2$ with the
product $\la v | w \ra_{L^2} := \int_{\Td} d x \, v  w
= (2 \pi)^{d} \sum_{k \in \Zd} \overline{v_k} w_k$ and
the corresponding norm $\|~ \|_{L^2}$.
For all  $n \in \reali$, the $n$-th Sobolev
space of zero mean functions on $\Td$ is
\parn
\vbox{
\beq \hz{n}(\Td) \equiv \hz{n} := \Big\{ v \in D'~|~
\langle \, v \rangle = 0~,
\sqrt{-\Delta}^{\,n} v \in L^2~ \}  \feq
$$ = \Big\{ v \in D'~|~v_0 = 0,~
\sum_{k \in \Zd \setminus \{0 \}} |k|^{2 n} |v_k|^2 < + \infty~\} $$
}
\noindent
(in the above $\langle v \rangle \in \reali$ is the mean
of $v$, i.e., by definition, the action
of $v$ on the test function $1/(2 \pi)^d$; moreover,
$\sqrt{-\Delta}^n v:= \sum_{k \in \Zd \setminus \{0 \}}
|k|^{n} v_k e_k$). $\hz{n}$ is a Hilbert space with the inner product
and the norm
\beq \la v | w \ra_n := \la \sqrt{-\Delta}^n v | \sqrt{-\Delta}^n w \ra_{L^2}
= (2 \pi)^d \! \! \! \! \! \sum_{k \in \Zd \setminus \{0 \}} \! \! \! |k|^{2 n} \overline{v_k} w_k ,
\quad \| v \|_n := \sqrt{\la v | v \ra_n}~; \feq
if $m \leqs n$, then $\hz{n} \subset \hz{m}$.
\salto
\textbf{The Laplacian and its semigroup.}
Let us consider the operator $\Delta := \sum_{rs=1}^d \partial_{s s} :
D' \vain D'$; of course $\Delta e_k = - |k|^2 e_k$ for all $k \in \Zd$.
For each $n \in \reali$, $\Delta$ maps continuously $\hz{n+2}$ into $\hz{n}$,
with $\| \Delta v \|_{n} \leqs \| v \|_{n+2}$ for all $v \in \hz{n+2}$.
We can define a semigroup $(e^{t \Delta})_{t \in [0,+\infty)}$ of linear operators on $D'$, putting
\beq e^{t \Delta} : D' \vain D'~, \qquad v \mapsto
e^{t \Delta} v := \sum_{k \in \Zd} e^{-| k |^2 t} v_k e_k~. \label{etd} \feq
Let $n \in \reali$. The following holds:
\beq e^{t \Delta} \hz{n} \subset \hz{n}~, \qquad \| e^{t \Delta} v \|_n \leqs e^{-t} \| v \|_n
\qquad \mbox{for $t \in [0,+\infty)$, $v \in \hz{n}$}~; \label{24} \feq
\beq e^{t \Delta} \hz{n-1} \subset \hz{n}~; \quad \exists \mu \in L^1((0,+\infty),
\reali)~ \mbox{(independent of $n$) such that}  \label{eqmu} \feq
$$ \| e^{t \Delta} v \|_n \leqs
\mu(t)  \| v \|_{n-1}~\mbox{for $t \in (0,+\infty)$, $v \in \hz{n-1}$} $$
(for the proof of \rref{eqmu} see, e.g., \cite{accau}, that
also gives an explicit expression for $\mu$ implying $\mu(t) = O(1/\sqrt{t})$
for $t \vain 0^{+}$ and $\mu(t) = e^{-t}$ for $t$ large).
The map $(t, v) \mapsto e^{t \Delta} v$ is continuous from $[0,+\infty) \times \hz{n}$ to
$\hz{n}$ and from $(0,+\infty) \times \hz{n-1}$ to
$\hz{n}$.
Moreover,
\beq v \in \hz{n+2} \Rightarrow (t \mapsto e^{t \Delta} v) \in  C([0,+\infty), \hz{n+2})
\cap C^1([0,+\infty), \hz{n})~, \label{delap} \feq
$$ {d \over d t} \big( e^{t \Delta } v \big) = \Delta (e^{t \Delta} v)
\qquad \mbox{for $t \in [0,+\infty)$}. $$
To go on, let
\beq
f \in C([0,+\infty), \hz{n+1})~; \feq
for each $t \in [0,+\infty)$ the function
$s \in (0,t) \mapsto e^{(t-s) \Delta} f(s)$
is in $L^1((0,t), \hz{n+2}$) because, on the grounds of \rref{eqmu},
$\| e^{(t-s) \Delta} f(s) \|_{n+2} \leqs \mu(t-s) \| f(s) \|_{n+1}$; therefore, the definition
\beq F(t) := \int_{0}^t d s \, e^{(t-s) \Delta} f(s) \qquad
\mbox{for $t \in [0,+\infty)$}
\label{26} \feq
produces a function $F \in C([0,+\infty), \hz{n+2})$. This function is also in  $C^1([0,+\infty), \hz{n})$
and satisfies an inhomogeneous heat equation with source term $f$:
\beq {d F \over d t}(t)
= \Delta F(t) + f(t)
\qquad \mbox{for $t \in [0,+\infty)$}~. \label{27} \feq
(Moreover, $F$ is the unique solution of \rref{27} with the above indicated regularity
and with $F(0)=0$.)
\salto
\textbf{Vector fields on $\Td$.}
Here and in the sequel,  ``a vector field on $\Td$'' means
``an $\reali^d$-valued distribution on $\Td$''. We
write $\Dv(\Td) \equiv \Dv$ for the space of such distributions;
these can be identified with $d$-tuples $v = (v_1,...,v_d)$, where
$v_r \in D'$ for each $r$.
Partial derivatives,
the Laplacian $\Delta$ and the operators $\sqrt{-\Delta}^n$,
$e^{t \Delta}$ are defined componentwise as maps from $\Dv$ to $\Dv$.
Any $v
\in \Dv$ has a weakly convergent Fourier expansion $v = \sum_{k \in \Zd} v_k e_k$, with
coefficients $v_k \in \complessi^d$ (such that $\overline{v_k} = v_{-k}$). \par
In the sequel $\mathbb{L}^p(\Td) \equiv \mathbb{L}^p$ denotes the space
of $L^p$ vector fields $\Td \vain \reali^d$.
For each  $n \in \reali$, the $n$-th Sobolev
space of zero mean vector fields on $\Td$ is
\par
\vbox{
\beq \Hz{n}(\Td) \equiv \Hz{n} := \Big\{ v \in \Dv~|~
\langle v \rangle = 0,~
\sqrt{-\Delta}^{\,n} v \in \mathbb{L}^2~ \}  \feq
$$ = \Big\{ v \in \Dv~|~v_0 = 0~,
\sum_{k \in \Zd \setminus \{0 \}} |k|^{2 n} |v_k|^2 < + \infty~\} $$
}
\noindent
(in the above, the mean $\langle v \rangle \in \reali^d$
is defined componentwise; we note that $\Hz{n} = \{ v = (v_1,...,v_d)
\in \Dv~|~v_r \in \hz{n}~\mbox{for each $r$} \}$).
$\Hz{n}$ is a Hilbert space with the inner product
and the norm
\beq \la v | w \ra_n := \la \sqrt{-\Delta}^n v | \sqrt{-\Delta}^n w \ra_{L^2}
= (2 \pi)^d \! \! \! \! \sum_{k \in \Zd \setminus \{0 \}} \! \! \! \! \! |k|^{2 n} \overline{v_k} \sc w_k,
~~ \| v \|_n := \sqrt{\la v | v \ra_n}~. \feq
Eqs. \rref{24}-\rref{27}
and the related statements have obvious analogues, where
$\hz{n}$ is replaced by $\Hz{n}$ for any $n$.
\salto
\textbf{Divergence free vector fields; the Leray projection.}
The space of divergence free vector fields on $\Td$ is
\beq \Ds := \{ v \in \bb{\DD}'~|~\dive \,v = 0 \}
= \{ v \in \bb{\DD}'~|~k \sc \, v_k = 0~\forall k \in \Zd~\}~. \feq
The \emph{Leray projection} is the linear, surjective map
\beq \LP : \Dv \vain \Ds~, \qquad
v \mapsto \LP v := \sum_{k \in \Zd} (\LP_k v_k) e_k~; \label{lp1} \feq
here $\LP_k$ is the orthogonal projection of $\complessi^d$ onto the orthogonal
complement of $k$, i.e.,
\beq \LP_\zz c = c~, \qquad \LP_k c = c - {k \sc \, c  \over | k |^2}\, k \quad \mbox{for $c \in \complessi^d$,
$k \in \Zd \setminus \{0 \}$}~. \label{lp2} \feq
For each $n \in \reali$, the \emph{$n$-th Sobolev space of zero mean,
divergence free vector fields} is
\beq \HM{n} := \Ds \cap \Hz{n}~; \label{hmn} \feq
this is a closed subspace of $\Hz{n}$, and thus
becomes a Hilbert space with the restriction of the inner product
$\la~|~\ra_n$.
One has
\beq \LP \Hz{n} = \HM{n}~, \quad \| \LP v \|_n \leqs \| v \|_n
\qquad \mbox{for $n \in \reali$, $v \in \Hz{n}$}~. \feq
The spaces \rref{hmn} are the basis of our treatment
of the NS equations; again, we have analogues
of Eqs. \rref{24}-\rref{27}
and of the related statements, where
$\hz{n}$ is replaced by $\HM{n}$ for any $n$.
\salto
\textbf{The NS bilinear map.}
Consider two vector fields $v, w$ on $\Td$ such that
$v \in \mathbb{L}^2$ and $\partial_s w \in \mathbb{L}^2$ for $s=1,...,d$;
then we have a well defined vector field $v \sc \partial w \in \mathbb{L}^1$
of components $(v \sc \partial w)_r := \sum_{s=1}^d v_s \partial_s w_r$ (which has
mean zero if $\dive \, v = 0$);
we can apply to this the Leray projection $\LP$ and form the (divergence free) vector field
\beq
\PPP(v, w) := - \LP(v \sc \partial w)~.
\feq
The bilinear map $\PPP$: $(v, w) \mapsto \PPP(v,w)$, which is
a main character of the incompressible NS equations, is known to
possess the following properties.
\begin{itemize}
\item[(i)]
For each $n > d/2$, $\PPP$ is continuous from
$\HM{n} \times \HM{n+1}$ to $\HM{n}$; so, there is a
constant $K_{n d} \equiv K_n$ such that
\beq \| \PPP(v, w) \|_n \leqs K_n \| v \|_n \| w \|_{n+1}
\quad \mbox{for $v \in \HM{n}$, $w \in \HM{n+1}$}~. \label{basic} \feq
\item[(ii)]
For each $n > d/2+1$, there is a constant $G_{n d} \equiv G_n$
such that
\beq |\la \PPP(v,w) | w \ra_n | \leqs G_n \| v \|_n
\| w \|^{2}_{n} \quad \mbox{for $v \in \HM{n}$, $w \in \HM{n+1}$}
\label{katineq} \feq
(this result is due to Kato, see \cite{Kato}).
\end{itemize}
From here to the end of the paper,
$K_n$ and $G_n$ are constants fulfilling the previous inequalities
(and not necessarily sharp). From \cite{cog} \cite{cok} we know that
we can take
\beq K_3 = 0.323~, \qquad G_3 = 0.438 \qquad \mbox{if $d=3$}~; \label{k3g3} \feq
these values will be useful in the sequel.
\salto
\textbf{The NS Cauchy problem.} From here to the end of the paper,
we fix a Sobolev order
\beq n \in \big({d \over 2} + 1, + \infty\big)~. \label{propen} \feq
Let us choose a Reynolds number $R \in [0,+\infty)$ and an initial datum
\beq \uz \in \HM{n+2}~. \label{in}\feq
\begin{prop}
\textbf{Definition.}
\textsl{
The (incompressible) NS Cauchy problem
with Reynolds number $R$ and initial datum $\uz$ is the following:
\beq \mbox{Find}~
u \in C([0, T), \HM{n+2}) \cap C^1([0,T), \HM{n}) \quad \mbox{such that} \label{cau} \feq
$$ {d u \over d t} = \Delta u + R \, \PPP(u,u)~, \qquad u(0) = \uz $$
(with $T \in (0, + \infty]$, depending on $u$)}.
\end{prop}
It is known \cite{Kat2} that the above Cauchy problem has a unique
maximal (i.e., not extendable) solution; any solution
is a restriction of the maximal one.
({\footnote{It is known as well that, if $u$ solves \rref{cau},
the function $(x,t) \mapsto u(x,t)$ is smooth on $\Td \times (0, T)$
(see, e.g., Theorem 15.2 (A) of \cite{Lem}); one could give stronger regularity
results with suitable assumptions on $\uz$. In the sequel we will not be
interested in such regularity matters.}})
\salto
\textbf{Approximate solutions of the NS Cauchy problem.}
In this paragraph we briefly rephrase
some basic results of \cite{appeul} with the notations
of the present paper (note that in \cite{appeul}
the NS equations were written in the form \rref{eulnu} rather than \rref{eul}).
We consider again the Cauchy problem \rref{cau}, for
given $n, R, \uz$ as above;
the definitions and the theorem that follow are reported from
\cite{appeul}, with obvious adaptations.
\begin{prop}
\textbf{Definition.}
\textsl{
An \emph{approximate solution} of the problem \rref{cau} is
any map $\ua \in C([0, \Ta), \HM{n+2}) \cap C^1([0,\Ta), \HM{n})$
(with $\Ta \in (0,+\infty]$). Given such a function,
we stipulate (i) (ii)}. \par\noindent
\textsl{
(i) The \emph{differential
error} of $\ua$ is
\beq e(\ua) := {d \ua \over d t} - \Delta \ua - R \, \PPP(\ua,\ua)~
\in C([0,\Ta), \HM{n})~;  \label{differr} \feq
the \emph{datum error} is
\beq \ua(0) - \uz \in \HM{n+ 2}~. \feq
(ii) Let $m \in \reali, m \leqs n$. A \emph{differential error estimator} of order $m$ for $\ua$
is a function $\ep_m \in C([0,\Ta), [0,+\infty))$ such that
\beq \| e(\ua)(t) \|_m \leqs \ep_m(t)~\mbox{~~for $t \in [0,\Ta)$}~. \label{erest} \feq
Let $m \in \reali$, $m \leqs n + 2$.
A \emph{datum error estimator} of order $m$ for $\ua$ is a real number
$\delta_m \in [0,+\infty)$ such that
\beq \| \ua(0) - \uz \|_m \leqs \delta_m~; \label{uaz}
\feq
a \emph{growth estimator} of order $m$ for $\ua$ is a function $\Dd_m \in C([0,\Ta), [0,+\infty))$
such that
\beq \| \ua(t) \|_m \leqs \Dd_m(t)~\mbox{~~for $t \in [0,\Ta)$}~.
\label{din} \feq
In particular the function $\ep_m(t) := \| e(\ua)(t) \|_m$,
the number $\delta_m := \| \ua(0) - \uz \|_m$ and
the function $\Dd_m(t) := \| \ua(t) \|_m$ will be called the
\emph{tautological} estimators of order $m$ for the differential error,
the datum error and the growth of $\ua$.}
\end{prop}
Let us observe that Eq. \rref{differr} could be read as follows:
the function $\ua$ fulfills the NS equations with an external
forcing $e(\ua)$. The \textsl{a posteriori} analysis of Chernyshenko \textsl{et al.}
for the approximate NS solutions is based on this remark and
on a robustness result for the NS equations with respect
to external forces (see \cite{Che};
more detailed information on the relations between this paper and
our approach is given in \cite{appeul}). \par
From here to the end of the section we consider
an approximate solution $\ua$ of the problem \rref{cau},
with domain $[0, \Ta)$; this is assumed to possess
differential, datum error and growth estimators
of orders $n$ or $n+1$, indicated with $\ep_n, \delta_n, \Dd_n,
\Dd_{n+1}$. \par
\begin{prop}
\textbf{Definition.}
\textsl{
Let $\Rr_n \in C([0,\Tc), \reali)$,
with $\Tc \in (0,\Ta]$. This function is said to fulfill the
\emph{control inequalities} if
\parn
\vbox{
\beq {d^{+} \Rr_n \over d t} \geqs - \Rr_n
+ R (G_n \Dd_n + K_n \Dd_{n+1}) \Rr_n + R \, G_n \Rr^2_n + \ep_n
\quad \mbox{in $[0,\Tc)$}, \label{cont1} \feq
\beq \Rr_n(0) \geqs \delta_n~. \label{cont2} \feq
}
In the above $d^{+}/ d t$ indicates the right, upper Dini derivative: so,
for all $t \in [0,\Tc)$,
$(d^{+} \Rr_n/ d t)(t) := \limsup_{h \vain 0^{+}} [\Rr_n(t+h) - \Rr_n(t)]/h$.}
\end{prop}
One can show that any function
$\Rr_n \in C([0, \Tc), \reali)$ fulfilling
the control inequalities
is automatically nonnegative
({\footnote{In fact the zero function
fulfills relations analogous to
\rref{cont1} \rref{cont2},
with $\geqs$ replaced by $\leqs$; by
standard comparison results, this implies $\Rr_n(t)
\geqs 0$ for all $t \in [0,\Tc)$.
The comparison results required to prove the last
statement have been reviewed in
\cite{appeul} (see Lemma 4.3 in the cited paper
and the related references).}}).
\begin{prop}
\label{main}
\textbf{Proposition.}
\textsl{
Assume
there is a function $\Rr_n \in C([0,\Tc), [0,+\infty))$ fulfilling
the control inequalities; consider the maximal solution $u$ of the
NS Cauchy problem \rref{cau}, and denote its domain with $[0,T)$.
Then
\beq T \geqs \Tc~, \label{tta} \feq
\beq \| u(t) - \ua(t) \|_n \leqs \Rr_n(t) \qquad \mbox{for $t \in [0,\Tc)$}~.
\label{furth} \feq
}
\end{prop}
\textbf{Proof (sketch).} For the sake of brevity, we put
\beq w(t) := u(t) - \ua(t) \qquad \mbox{ for $t \in [0, \min(T, \Ta))$}. \feq
The function $w$ fulfills
${d w/d t} = \Delta w + R \PPP(\ua, w) + R \PPP(w, \ua) + R \PPP(w,w) - e(\ua)$ and
$w(0) = \uz - \ua(0)$.
From here, using the inequalities \rref{basic} \rref{katineq}
about $\PPP$ and \rref{erest} \rref{uaz} \rref{din} about $e(\ua)$, $\uz - \ua(0)$ and
$\ua$ one infers that the continuous function $\| w \|_n : t \in [0, \min(T, \Ta)) \mapsto
\| w(t) \|_n$ fulfills
\beq {d^{+} \| w \|_n \over d t} \leqs - \, \| w \|_n + R (G_n \Dd_n + K_n \Dd_{n+1}) \| w \|_n + R
G_n \| w \|^2_n + \ep_n ~, \label{ver1} \feq
\beq \| w(0) \|_n \leqs \delta_n \label{ver2} \feq
(note that the usual derivative $d \| w \|_n/d t$ might fail to exist at times $t$ such that $w(t)=0$;
this is the reason for using the Dini derivative). From \rref{ver1} \rref{ver2},
from the control inequalities \rref{cont1} \rref{cont2} and from
a comparison theorem ({\footnote{See again Lemma 4.3 of \cite{appeul}.}}) one infers
\beq \| w(t) \|_n \leqs \Rr_n(t) \qquad \mbox{for $t \in [0, \min(T, \Tc))$} \feq
Finally, one has
\beq \min(T, \Tc) = \Tc~; \feq
in fact, if $T <  \Tc$ we would have
$\| u(t) \|_n \leqs \| w(t) \|_n + \| \ua(t) \|_n \leqs \Rr_n(t) + \Dd_n(t)$
for all $t \in [0,T)$ and this would imply $\limsup_{t \to T^{-}} \| u(t) \|_n
\leqs \Rr_n(T) + \Dd_n(T) < +\infty$, contradicting standard results on the
maximal solution of the NS Cauchy problem.
\fine
Paper \cite{appeul} presents some applications of
the previous proposition, where $\ua$ is constructed
by the Galerkin method. (For completeness
we mention that the framework of \cite{appeul}
also covers the case of the Euler equations,
i.e., the limit case $\nu \vain 0$ of \rref{eulnu}
which is formally equivalent to $R  \vain +\infty$; some
applications to the Euler equations have been considered both
in \cite{appeul} and in \cite{padova}.) \par
In the sequel we present an application of Proposition \ref{main},
choosing for $\ua$ a polynomial in $R$
(see Eq.\rref{tay}).
In the next two sections we develop this approach
in general terms, giving the error estimators
for an approximate solution of this kind;
in the final section we apply this procedure
choosing for $\uz$
the BNW initial datum.
\section{Reynolds number expansions
as approximate NS solutions}
\label{reynolds}
Let us recall that $n \in (d/2+1,+\infty)$, and consider
the NS Cauchy problem \rref{cau} with $R \in [0,+\infty)$
and datum $\uz \in \HM{n+2}$.
Let us choose an order $N \in \{0,1,2,...\}$ and consider as an approximate
solution for \rref{cau} a polynomial of degree $N$ in $R$, of the form
\parn
\vbox{
\beq u^N : [0, +\infty) \vain \HM{n+2}~,
\qquad t \mapsto u^N(t) := \sum_{j=0}^N R^j u_j(t)~, \label{defun} \feq
$$ u_j \in C([0,+\infty), \HM{n+2})
\cap C^1([0,+\infty), \HM{n}) \qquad \mbox{for $j=0,...,N$}; $$
}
the functions $u_j$ herein are to be determined.
\begin{prop}
\label{un}
\textbf{Proposition.}
(i) Let $u^N$ be as in \rref{defun}. The datum and differential errors of $u^N$ are
\beq u^N(0) - \uz = (u_0(0) - \uz) + \sum_{j=1}^N R^j u_j(0)~; \label{daer} \feq
\beq e(u^N) = \Big({d u_0 \over d t} - \Delta u_0\Big) +
\sum_{j=1}^{N} R^j \Big[ {d u^j \over d t} - \Delta u_j - \sum_{\ell=0}^{j-1} \PPP(u_\ell, u_{j-\ell-1}) \Big]
\label{eun} \feq
$$ - \sum_{j=N+1}^{2 N+1} R^j \! \! \! \! \! \sum_{\ell=j-N-1}^N \PPP(u_\ell, u_{j - \ell-1}). $$
(ii) One can define recursively a
family of functions $u_j \in C([0,+\infty), \HM{n+2})$ $\cap$ $C^1([0,+\infty), \HM{n})$
setting
\beq u_0(t) := e^{t \Delta} \uz \qquad \mbox{for $t \in [0,+\infty)$}~, \label{recurz} \feq
\beq u_{j}(t) := \sum_{\ell=0}^{j-1}
\int_{0}^t d s \, e^{(t-s) \Delta} \PPP(u_{\ell}(s), u_{j-\ell-1}(s))
\quad \mbox{for $t \in [0,+\infty)$, $j=1,...,N$}. \label{recur} \feq
With this choice of $u_0,...,u_N$ the coefficients of $R^0, R^1,..., R^N$ in Eqs. \rref{daer}
and \rref{eun} vanish, so that
\beq u^N(0) - \uz = 0~; \label{daerr} \feq
\beq e(u^N)
= - \sum_{j=N+1}^{2 N+1} R^j  \sum_{\ell=j-N-1}^N \PPP(u_\ell, u_{j - \ell-1})~.
\label{eunn} \feq
The second equation implies
\par
\vbox{
\beq \| e(u^N)(t) \|_n
\leqs K_n \sum_{j=N +1}^{2 N +1} R^j
\sum_{\ell=j-N -1}^N \| u_\ell(t) \|_n \| u_{j - \ell -1}(t) \|_{n+1}
~\mbox{for $t \in [0,+\infty)$}~.
\label{est} \feq
}
\end{prop}
\textbf{Proof.}
(i) Eq. \rref{daer} is obvious. Let us prove Eq. \rref{eun};
to this purpose, we note that
$$ e(u^N) = {d u^N \over d t} - \Delta u^N - R \, \PPP(u^N, u^N)  =
({d \over d t} - \Delta)\big( \sum_{j=0}^N R^j u_j \big)
- R \, \PPP\big(\sum_{\ell=0}^N R^\ell u_\ell,
\sum_{h=0}^N R^h u_h\big) $$
$$ = \sum_{j=0}^N R^j ({d u_j \over d t} - \Delta u_j)
- \sum_{\ell,h=0}^N R^{\ell+h+1} \PPP(u_\ell, u_h) $$
$$ = ({d u_0 \over d t} - \Delta u_0) +
\sum_{j=1}^N R^j ({d u_j \over d t} - \Delta u_j) - \sum_{j=1}^{2 N+1}
R^j \! \! \! \sum_{(\ell,h) \in I_{N j}} \PPP(u_\ell, u_h)~, $$
$$ I_{N j} := \{ (\ell,h) \in \{0,...,N\}^2~|~\ell+h+1= j \}~. $$
One easily checks that
$$ j \in \{1,...,N\}~\Rightarrow~
I_{N j} = \{ (\ell,j-\ell-1)~|~\ell \in \{0,...,j-1\} \}~, $$
$$ j \in \{N+1,...,2 N+1\}~\Rightarrow~
I_{N j} = \{ (\ell,j-\ell-1)~|~\ell \in \{j-N-1,...,N\} \}~; $$
this readily yields the thesis \rref{eun}. \parn
(ii) First of all, let us prove by recurrence over $j$ that the functions
$u_j$ of Eqs. \rref{recurz} \rref{recur} are well defined and belong to $C([0,+\infty), \HM{n+2})
\cap C^1([0,+\infty), \HM{n})$, for all $j \in \{0,...,N \}$. \parn
For $j=0$, this follows from an obvious vector analogue of
the regularity statement in \rref{delap}. Now, let $j \in \{1,...,N\}$ and
assume the thesis to hold up to the order $j-1$; then
the functions
$t \in [0,+\infty) \mapsto \PPP(u_{\ell}(t), u_{j-\ell-1}(t))$ ($\ell = 0,....,j-1$)
are in $C([0,+\infty), \HM{n+1})$ due to the properties of $\PPP$.
By obvious vector analogues of the considerations accompanying Eqs.
\rref{26} \rref{27}, with $f$ replaced by $\PPP(u_{\ell}, u_{j-\ell-1}) $, we see that
the functions $t \mapsto \int_{0}^t d s \, e^{(t-s) \Delta} \PPP(u_{\ell}(s), u_{j-\ell-1}(s))$
belong to $C([0,+\infty) \HM{n+2}) \cap C^1([0,+\infty), \HM{n})$; the same can be said of $u_{j}$, which is
the sum over $\ell$ of these functions. \par
To go on, let us note that the definition \rref{recurz}
implies
\beq u_0(0) = \uz~, \qquad {d u_0 \over d t} = \Delta u_0~, \label{calore} \feq
while the definition \rref{recur} implies
\beq u_j(0) = 0~,~
{d u_{j} \over d t}(t) = \Delta u_{j}(t) + \sum_{\ell=0}^{j-1} \PPP(u_{\ell}(t), u_{j-\ell-1}(t))
~~\mbox{for}~ j \in \{1,...,N \} \label{calorej} \feq
(to check the last equality, use Eq. \rref{27} with $f$ replaced by
$\sum_{\ell} \PPP(u_{\ell}, u_{j-\ell-1})$). Eqs. \rref{calore}
and \rref{calorej} indicate, respectively, the vanishing
of the coefficients of $R^0$ and $R^j$ ($j=1,...,N$) in both
Eqs. \rref{daer} \rref{eun}. \par
Now Eqs. \rref{daerr} \rref{eunn} are obvious.
Eq. \rref{eunn} and the inequality
$\| \PPP(u_\ell(t), u_{j - \ell-1}(t)) \|_n$
$\leqs K_n \| u_\ell(t) \|_n \| u_{j - \ell -1}(t) \|_{n+1}$
gives immediately Eq. \rref{est}. \fine
\begin{rema}
\textbf{Remark.} Apart from technicalities, the previous result
can be described as follows. The coefficients of $R^0, R^1,..., R^N$
in the expressions \rref{daer} \rref{eun} for the datum and differential errors of $u^N$
are zero if $u_0$ fulfills the heat equation with initial datum
$\uz$ and $u_j$ fulfills, for each $j \in \{1,...,N\}$, an inhomogeneous
heat equation with zero initial datum and a source term depending on $u_{0}, ..., u_{j-1}$.
The positions \rref{recurz} \rref{recur} provide the solutions for these initial value problems.
\fine
\end{rema}
Let $N \in \{0,1,2,...\}$; we define again
$u^N$ and
$u_j$ ($j=0,...,N$) via Eqs. \rref{defun} \rref{recurz} \rref{recur},
and regard $u^N$ as an approximate solution of the Cauchy problem \rref{cau}, of
domain $[0,+\infty)$. This has
the tautological datum error and growth estimators
\beq \delta_n := 0, \qquad \Dd_n(t) := \| u^N(t) \|_n, \quad
\Dd_{n+1}(t) := \| u^{N}(t) \|_{n+1}~. \label{vi} \feq
The differential error $e(u^N)$ can be
expressed via Eq. \rref{eunn}; we can use for it
the tautological estimator
\beq \ep_n(t) := \| e(u^N)(t) \|_n
\label{ep1} \feq
or the rougher estimator indicated by \rref{est}, i.e., the function
\beq \ep_n(t) := K_n \sum_{j=N +1}^{2 N +1} R^j
\sum_{\ell=j-N -1}^N \| u_\ell(t) \|_n \| u_{j - \ell -1}(t) \|_{n+1}~. \label{ep2} \feq
Now, using Proposition \ref{main} with $\ua = u^N$
and the above estimators we obtain the following.
\begin{prop}
\textbf{Corollary.}
\textsl{Let $N \in \{0,1,2,...\}$; define
$\delta_n, \Dd_n, \Dd_{n+1}$ via Eq. \rref{vi} and
$\ep_n$ via Eq. \rref{ep1} or \rref{ep2}, for $t \in [0,+\infty)$.
Suppose there is a function
$\Rr_n \in C([0,\Tc), [0,+\infty))$,
with $\Tc \in (0,+\infty]$, fulfilling the
control inequalities \rref{cont1} \rref{cont2}.
Consider the maximal solution $u$ of the
NS Cauchy problem \rref{cau}, of domain $[0,T)$;
then
\beq T \geqs \Tc~, \qquad \| u(t) - u^N(t) \|_n \leqs \Rr_n(t) \qquad \mbox{for $t \in [0,\Tc)$}~.
\label{furthe} \feq
In particular, the control inequalities \rref{cont1} \rref{cont2}
are satisfied by any function $\Rr_n \in C^1([0,\Tc), [0,+\infty))$
fulfilling the \emph{control Cauchy problem}
\beq {d \Rr_n \over d t} = - \Rr_n
+ R (G_n \Dd_n + K_n \Dd_{n+1}) \Rr_n + R \, G_n \Rr^2_n + \ep_n
~\mbox{on $[0,\Tc)$},~~\Rr_n(0) =0. \label{concau} \feq
}
\end{prop}
\begin{rema}
\textbf{Remarks.} (i) Even though rougher than the
estimator \rref{ep1}, the estimator \rref{ep2} is interesting because
its computation is less expensive; this is a relevant fact,
especially in applications with a large $N$. \parn
(ii) Let us consider the estimator \rref{ep1}, writing via \rref{eunn}
the error therein; alternatively, let us use the estimator \rref{ep2}.
In both cases, it is natural to write
\beq \ep_n(t) = R^{N+1} \tilde{\ep}_n(t) \qquad \mbox{for $t \in [0,+\infty)$}~, \feq
where $\tilde{\ep}(t)$ is a suitable function;
this is a polynomial of degree $N$ in $R$ in the case \rref{ep2},
and the square root of a polynomial of degree $2 N$ in $R$ in the case \rref{ep1}.
Consequently, the solution of
the Cauchy problem \rref{concau} can be written as
\beq \Rr_n(t) = R^{N+1} \tilde{\Rr}_n(t)~, \feq
where $\tilde{\Rr}_n \in C^1([0,\Tc), [0,+\infty))$ is such that
\beq {d \tilde{\Rr}_n \over d t} = - \tilde{\Rr}_n
+ R (G_n \Dd_n + K_n \Dd_{n+1}) \tilde{\Rr}_n + R^{N+2} G_n \tilde{\Rr}^2_n + \tilde{\ep}_n
~,~~\tilde{\Rr}_n(0) =0. \label{concaur} \feq
(iii) Eq. \rref{furthe} has a number of obvious implications; let us give two examples. \parn
(iii$_1$) The inequalities $\| u^N(t) \|_n - \| u(t) - u^N(t) \|_n \leqs
\| u(t) \|_n \leqs \| u^N(t) \|_n + \| u(t) - u^N(t) \|_n$, the
definition of $\Dd_n(t)$ in \rref{vi} and Eq. \rref{furthe} for $\Rr_n$ imply
\beq \Dd_n(t) - \Rr_n(t) \leqs \| u(t) \|_n \leqs \Dd_n(t) + \Rr_n(t)
\quad \mbox{for $t \in [0, \Tc)$}~. \label{conc1} \feq
(iii$_2$) For any $v \in \Hz{n}$, the equation $\| v \|^{2}_n =
(2 \pi)^d \sum_{k \in \Zd \setminus \{0 \}} |k|^{2 n} | v_k|^2$ implies
$(2 \pi)^{d/2} | v_k |$ $\leqs \| v \|_n/|k|^n$ for all $k$. This inequality
with $v = u(t) - u^N(t)$ and \rref{furthe} give
\beq (2 \pi)^{d/2} | u_k(t) - u^N_k(t) | \leqs {\Rr_n(t) \over |k|^n}
\qquad \mbox{for $k \in \Zd \setminus \{0 \}$ and $t \in [0,\Tc)$}~. \feq
\end{rema}
\section{Implementing the recursion relations \rref{recurz} \rref{recur}}
\label{imple}
The recursion relations mentioned in the title are the main characters
of Proposition \ref{un}; their building blocks are, essentially, the linear maps $\UU, \KK$
defined as follows:
\beq \UU : \hz{n+2} \vain C([0,+\infty), \hz{n+2}) \cap C^1([0,+\infty), \hz{n}),
\qquad z \mapsto \UU z \label{uz} \feq
$$ \mbox{where}~~ \UU z : [0,+\infty) \vain \hz{n+2}~, \qquad t \mapsto (\UU z)(t) := e^{t \Delta} z~; $$
\beq \KK : C([0,+\infty), \hz{n+1}) \vain C([0,+\infty), \hz{n+2})
\cap C^1([0,+\infty), \hz{n})~, \qquad f \mapsto \KK f \label{kf} \feq
$$ \mbox{where}~~ \KK f : [0,+\infty) \vain \hz{n+2}~, \qquad t \mapsto (\KK f)(t) :=
\int_{0}^t d s \, e^{(t-s) \Delta} f(s) $$
(as for the domains and codomains of the above maps, recall Eqs. \rref{delap} \rref{27}
and the related comments). These maps have vector analogues, denoted
for simplicity with the same letters,
\beq \UU : \Hz{n+2} \vain C([0,+\infty), \Hz{n+2}) \cap C^1([0,+\infty),\Hz{n}),
\qquad z \mapsto \UU z~, \label{uzv}\feq
\beq \KK : C([0,+\infty), \Hz{n+1}) \vain C([0,+\infty), \Hz{n+2}) \cap C^1([0,+\infty),
\Hz{n})~, \qquad f \mapsto \KK f~; \label{kfv} \feq
these are defined rephrasing Eqs. \rref{uz} \rref{kf}.
In the
sequel, it will always be clear from the context whether
the symbols $\UU, \KK$ refer to the maps \rref{uz} \rref{kf}
or to the maps \rref{uzv} \rref{kfv}.
For $z = (z_1,...,z_d) \in \HM{n+2}$ and $f = (f_1,...,f_d)
\in C([0,+\infty), \Hz{n+1})$ we have
\beq \UU z = (\UU z_1, ..., \UU z_d)~, \qquad
\KK f = (\KK f_1,..., \KK f_d)~. \feq
With the above notations, for any $\uz \in \HM{n+2}$
the recursion relations \rref{recurz} \rref{recur} can be written as
\beq u_0 := \UU \uz~, \label{recurrez} \feq
\beq u_{j} := \sum_{\ell=0}^{j-1} \KK \PPP(u_{\ell}, u_{j-\ell-1})
\label{recurre} \qquad (j=1,...,N)~. \feq
Hereafter we give a set of elementary computational rules
about the maps $\UU, \KK$ and $\PPP$;
these allow a straightforward computation of the functions
$u_0, u_1,...,u_N$ in Eqs. \rref{recurrez} \rref{recurre},
at least in the case when $\uz$ is a Fourier polynomial
(see the remark at the end of this section).
Our computational rules involve the functions
$e_k \in D$, $e_k(x) := e^{i k \sc x}$ ($k \in \Zd$) and
\beq B_{a, b} : [0,+\infty) \vain \reali~,
\qquad t \mapsto B_{a, b}(t) := t^a e^{-b t} \qquad (a,b \in \naturali)~, \feq
as well as the products
\beq B_{a, b} \, e_k : [0,+\infty) \vain C^{\infty}(\Td, \complessi)~, \qquad t \mapsto t^a e^{-b t} \, e_k~.\feq
The formulation of our rules entails the complexified versions of the maps
$\UU$, $\KK$ and $\PPP$, denoted with the same symbols. For example, the complexified versions
of the maps $\UU$ and $\KK$ in Eqs. \rref{uz} \rref{kf} are defined as in the cited equations replacing
systematically the spaces $D'$ and $H^n_{0}$ with their complex
analogues $D'_{\complessi}$ (the space of all complex distributions
on $\Td$) and $H^{n}_{0 \complessi}$ (the space of zero mean complex distributions
$v$ such that $\sqrt{-\Delta}^{n} \, v$ is square integrable).
With the previous warnings, we summarize our computational rules in the forthcoming
lemma.
\begin{prop}
\textbf{Lemma.} Let $a, a',b, b' \in \naturali$, $k, k' \in \Zd$ and $s \in \{1,.., d \}$.
Then
\beq \partial_s (B_{a, b} \, e_k) = i k_s B_{a, b} \, e_k~; \label{410} \feq
\beq (B_{a, b} \, e_k) (B_{a', b'} \, e_{k'}) = B_{a + a', b + b'} \, e_{k + k'}~.
\label{411} \feq
In addition, let $k \neq 0$ (so that the functions below are in the
domains of the complexified maps $\UU$ or $\KK$). Then
\beq \UU e_k = B_{0, |k|^2} \, e_k~; \label{412} \feq
\beq {~} \hspace{-0.4cm} \KK (B_{a, b} \, e_k) = a! \Big( {B_{0, |k|^2} \over (b - |k|^2)^{a+1}}
- \sum_{\ell=0}^a {B_{\ell, b} \over  (b - |k|^2)^{a + 1 - \ell} \ell!} \Big) \, e_k
~~~\mbox{if $b \neq |k|^2$}; \label{413} \feq
\beq \KK(B_{a, |k|^2} \, e_k) = {B_{a +1, |k|^2} \over a + 1}  \, e_k~. \label{414}\feq
\end{prop}
\textbf{Proof.} Eqs. \rref{410} \rref{411} are obvious consequences of the
definitions of $e_k$ and $B_{a, b}$. Eq. \rref{412} is just a reformulation of
the relation $e^{t \Delta} e_k = e^{- |k|^2 t} \, e_k$. \par
In order to derive Eqs. \rref{413} \rref{414}, we note that
(both for $b \neq |k|^2$ and for $b = |k|^2$),
\parn
\vbox{
\beq \KK (B_{a, b} \, e_k)(t) = \int_{0}^t d s \, e^{(t-s) \Delta} (s^a e^{-b s} \, e_k)  \label{kb} \feq
$$ = \int_{0}^t d s \, s^a e^{-b s} e^{-|k|^2 (t-s)} \, e_k =
e^{-|k|^2 t} \left(  \int_{0}^t d s \, s^a e^{-(b - |k|^2) s} \right) \, e_k~. $$
}
\parn
\vbox{
\noindent
Let $b \neq |k|^2$. Then, a change of variables $\sigma = (b - |k|^2) s$ gives
$\int_{0}^t d s \, s^a e^{-(b - |k|^2) s}$ $= 1/(b - |k|^2)^{a+1} \int_{0}^{(b - |k|^2) t} d \sigma \, \sigma^a e^{-\sigma}$;
thus
}
\beq \KK (B_{a, b} \, e_k)(t) = {e^{-|k|^2 t} \over (b - |k|^2)^{a+1}} \, \gamma(a+1, (b - |k|^2) t) \, e_k~,
\label{kg} \feq
where we have introduced the incomplete Gamma function
\beq \gamma(\alpha, y) := \int_{0}^y
d \sigma \, \sigma^{\alpha-1} e^{-\sigma} \qquad \mbox{for $y \in \reali$, $\alpha \in \{1,2,3,...\}$}~.
\feq
It is known that
\beq \gamma(a+1, y) = a! \Big( 1 - e^{-y} \sum_{\ell=0}^a {y^\ell \over \ell!} \Big) \feq
(see, e.g., \cite{Nist}, page 177, Eqs. (8.4.7) and (8.4.11)). Inserting
this result into \rref{kg} we obtain that, for $b \neq |k|^2$,
\beq \KK (B_{a, b} \, e_k)(t) =
a! \left(
{e^{- |k|^2 t} \over (b - |k|^2)^{a+1}}-
\sum_{\ell=0}^a {t^{\ell} e^{-b t} \over (b - |k|^2)^{a+1 - \ell} \ell!} \,
\right) \, e_k~; \feq
this proves Eq. \rref{413}.
Finally, in the case $b = |k|^2$ we obtain from \rref{kb} that
\beq \KK (B_{a, |k|^2} \, e_k)(t) = e^{-|k|^2 t}
\left(  \int_{0}^t d s \, s^a  \right) \, e_k = {t^{a+1} e^{-|k|^2 t} \over a + 1} \, e_k~; \feq
this proves Eq. \rref{414}. \fine
\begin{rema}
\textbf{Remark.}
Let us return to the recursion rules \rref{recurrez} \rref{recurre},
assuming that the initial datum $\uz$ is a Fourier polynomial: by
this we mean that $\uz$ has finitely many non zero Fourier coefficients.
In this case, due to Eqs. \rref{recurrez} and \rref{412}
all components of $u_0$ are linear combinations
of finitely many functions of the form $B_{0, |k|^2} \, e_k$. Using Eq.
\rref{recurre} with the results of the previous lemma
we see that, for $j=1,...,N$, each component of $u_j$
is a finite linear combination of functions of the form $B_{a, b} \, e_k$.
We note that each term
$\PPP(u_{\ell}, u_{j-\ell-1})$ in Eq. \rref{recurre} can be
calculated as follows: first of all, one computes
each component of $u_{\ell} \sc \partial u_{j - \ell-1}$
using elementary considerations of bilinearity, together with
Eqs. \rref{410} \rref{411}; next, one obtains
$\PPP(u_{\ell}, u_{j-\ell-1}) = - \LP(u_{\ell} \sc \partial u_{j - \ell-1})$
using the expression \rref{lp1}\rref{lp2} for the Leray projection.
\end{rema}
\section{An application with the
Behr-Ne$\check{\mbox{c}}$as-Wu (BNW) initial datum}
\label{behr}
Throughout this section we work with
\beq d = 3 ~, \qquad n= 3 \feq
and any Reynolds number $R \in [0,+\infty)$. The Cauchy problem
\rref{cau} takes the form
\parn
\vbox{
\beq \mbox{\textsl{Find}}~
u \in C([0, T), \HM{5}) \cap C^1([0,T), \HM{3}) \quad \mbox{\textsl{such that}} \label{caun} \feq
$$ {d u \over d t} = \Delta u + R \, \PPP(u,u)~, \qquad u(0) = \uz~. $$
}
The initial datum $\uz$ in $\HM{5}$ (in fact, in $\HM{m}$ for any real $m$)
is chosen as follows:
\parn
\vbox{
\beq \uz(x_1, x_2, x_3) := 2 \big(\cos(x_1 + x_2) + \cos(x_1 + x_3), \label{unec0}\feq
$$ - \cos(x_1 + x_2) + \cos(x_2 + x_3),
 - \cos(x_1 + x_3) - \cos(x_2 + x_3)\big)~. $$
 }
This is a Fourier polynomial; indeed
\par
\vbox{
\beq \uz = \sum_{k = \pm a, \pm b, \pm c} u_{* \, k} e_k~, \label{unec} \feq
$$ a := (1,1,0),~~ b := (1,0,1),~~ c := (0,1,1)~; $$
$$ u_{*, \pm a} := (1,-1,0)~, \quad u_{*, \pm b}
:= (1,0,-1)~, \quad u_{*, \pm c} := (0,1,-1)~. $$
}
\noindent
This initial datum has been introduced by Behr, Ne$\check{\mbox{c}}$as
and Wu in \cite{Nec}; these authors have considered the datum \rref{unec0} \rref{unec}
as the origin of a possible blow-up for the Euler equations
(i.e., for the zero viscosity limit of \rref{eulnu}). \par
We have disputed the BNW blow-up conjecture in \cite{bnw}; in the present paper, independently
of any opinion on the validity of this conjecture, we consider
the datum \rref{unec0} \rref{unec} in presence of viscosity and analyze it by the
Reynolds expansion  method of sections \ref{reynolds} and \ref{imple}.
In the final subsection we compare the outcomes of this approach
with the results on the BNW datum obtained for the viscous case
in \cite{appeul}, where a Galerkin approximate solution was
employed in place of the Reynolds expansion.
\salto
\textbf{Setting up the Reynolds expansion.} The expansion of order $N$ relies
on the function $u^N := \sum_{j=0}^N R^j u_j : [0,+\infty) \vain \HM{5}$
where $u_0,...,u_N$ are computed via the recursion rules
\rref{recurrez} \rref{recurre}, starting from the BNW datum \rref{unec};
after finding $u^N$, one computes
the related growth and error estimators $\Dd_3, \Dd_4, \ep_3$ and
sets up the control Cauchy problem
\beq {d \Rr_3 \over d t} = - \Rr_3
+ R (G_3 \Dd_3 + K_3 \Dd_{4}) \Rr_3 + R \, G_3 \Rr^2_3 + \ep_3
~\mbox{on $[0,\Tc)$},~~\Rr_3(0) =0~, \label{conca} \feq
with $K_3$ and $G_3$ as in \rref{k3g3}. In the sequel
\beq \Rr_3 : [0, \Tc) \vain [0,+\infty) \label{r3} \feq
always denotes the maximal solution of this problem
({\footnote{More precisely: one considers the real valued $C^1$ functions
fulfilling \rref{conca}, which are automatically
nonnegative (see the comment immediately after
Proposition \ref{main}); the maximal $C^1$ solution
(i.e., the solution with the largest domain) is the
function in \rref{r3}.}}). \par
The above conceptual scheme has been implemented
in the cases $N=1,2,5$, on which we report in the sequel.
The case $N=1$ is useful to illustrate
in full detail the method; the other ones
give more accurate results, at the price of
more expensive computations.
In all cases, for the practical computation of
$u_0,u_1,...$ and of the estimators $\Dd_3, \Dd_4, \ep_3$,
we have employed MATHEMATICA in the symbolic mode. Then
the control Cauchy problem has been solved numerically (using
again MATHEMATICA) for several sample values of $R$;
these numerical calculations are very realiable, since
they concern a simple one-dimensional ODE.
The outcomes of such computations give evidence for the
following picture.
\begin{itemize}
\item[(i)] \label{pagitem} There is a critical Reynolds number $R_{*}$ (depending on $N$) such that
$\Tc = +\infty$ for $0 \leqs R \leqs R_{*}$, and $\Tc < + \infty$ for $R > R_{*}$.
Moreover, for $0 \leqs R \leqs R_{*}$ one has $\Rr_3(t) \vain 0^{+}$ for
$t \vain +\infty$, while for $R > R_{*}$ one has $\Rr_3(t) \vain + \infty$
for $t \vain \Tc^{-}$.
\end{itemize}
On the grounds of our general setting, the results (i) on the control
Cauchy problem yield the following conclusion.
\begin{itemize}
\item[(ii)] Let us consider the maximal solution $u$ of the NS Cauchy problem \rref{caun}
and the Reynolds expansion (for $N=1,2$ or $5$), with its critical number $R_{*}$.
For $0 \leqs R \leqs R_{*}$, $u$ is global and $\| u(t) - u^N(t) \|_3
\leqs \Rr_3(t)$ for all $t \in [0,+\infty)$. For $R > R_{*}$, we
can only grant that the domain of $u$ contains the interval $[0, \Tc)$, and
that $\| u(t) - u^N(t) \|_3 \leqs \Rr_3(t)$ for all $t \in [0,\Tc)$.
\end{itemize}
In the sequel we give more detailed information about
each one of the cases $N=1,2,5$. As expected, when $N$ is
increased the critical value $R_{*}$ increases as well.
Our tests also indicate that, for a given $R$,
when $N$ increases the same happens of $\Tc$; on the contrary
$\Rr_3$ becomes smaller, thus giving
a more stringent estimate on $\| u(t) - u^N(t) \|_3$.
\salto
\textbf{Case $\boma{N=1}$.} This relies on
\beq u^1 := u_0 + R u_1 : [0,+\infty) \vain \HM{5}~. \feq
The expressions of $u_j = (u_{j, 1}, u_{j, 2}, u_{j, 3})$ for $j=0,1$
are as follows (with $B_{a, b}(t) := t^a e^{-b t}$,
as in the previous section):
\parn
\vbox{
\beq u_{0, 1} =
B_{0, 2}( e_{(-1, -1, 0)} + e_{(-1, 0, -1)} + e_{(1, 0, 1)} + e_{(1, 1, 0)})~, \label{u0} \feq
$$ u_{0, 2}  = B_{0, 2} (- e_{(-1, -1, 0)} + e_{(0, -1, -1)} +
e_{(0, 1, 1)} - e_{(1, 1, 0)})~,$$
$$ u_{0, 3}  = B_{0, 2} (- e_{(-1, 0, -1)} - e_{(0, -1, -1)} - e_{(0, 1, 1)} - e_{(1, 0, 1)})~; $$
}
\parn
\vbox{
\beq u_{1, 1}
= {2 i \over 3} B_{0, 4} (e_{(-2, -1, -1)} + e_{(-1, -2, -1)} - e_{(-1, -1, -2)} + e_{(1, 1, 2)} -
        e_{(1, 2, 1)} - e_{(2, 1, 1)}) \label{u1} \feq
$$ + {2 i \over 3} B_{0, 6} (-e_{(-2, -1, -1)} - e_{(-1, -2, -1)} + e_{(-1, -1, -2)} -
        e_{(1, 1, 2)} + e_{(1, 2, 1)} + e_{(2, 1, 1)})~, $$
$$ u_{1, 2} =
{2 i \over 3} B_{0, 4} (-e_{(-2, -1, -1)} - e_{(-1, -2, -1)} - e_{(-1, -1, -2)} +
        e_{(1, 1, 2)} + e_{(1, 2, 1)} + e_{(2, 1, 1)}) $$
$$ + {2 i \over 3} B_{0, 6} (e_{(-2, -1, -1)} + e_{(-1, -2, -1)} + e_{(-1, -1, -2)} - e_{(1, 1, 2)} -
        e_{(1, 2, 1)} - e_{(2, 1, 1)})~, $$
$$ u_{1, 3} =
{2 i \over 3} B_{0, 4} (-e_{(-2, -1, -1)} + e_{(-1, -2, -1)} + e_{(-1, -1, -2)} -
        e_{(1, 1, 2)} - e_{(1, 2, 1)} + e_{(2, 1, 1)}) $$
$$ + {2 i \over 3} B_{0, 6} (e_{(-2, -1, -1)} - e_{(-1, -2, -1)} - e_{(-1, -1, -2)} + e_{(1, 1, 2)} +
        e_{(1, 2, 1)} - e_{(2, 1, 1)})~. $$
}
The next step is to compute the differential error $e(u^1)$
(which can be expressed via \rref{eunn}) and the (time dependent) norms
\beq \Dd_3 := \| u^1 \|_3, \quad
\Dd_{4} := \| u^1 \|_{4}, \quad
\ep_3 := \| e(u^1)\|_3~. \feq
These are as follows:
\beq \Dd_3 = 4 \sqrt{6} \, (2 \pi)^{3/2}
\Big[B_{0, 4} + 18 R^2 (B_{0, 8} - 2 B_{0, 10} + B_{0, 12}) \Big]^{1/2} \feq
\beq \Dd_4 = 8 \sqrt{3} \, (2 \pi)^{3/2}
\Big[B_{0, 4} + 54 R^2 (B_{0, 8} - 2 B_{0, 10} + B_{0, 12}) \Big]^{1/2}~, \feq
\parn
\vbox{
\beq \ep_3 = 8 \sqrt{2/3} \, (2 \pi)^{3/2} R^2 \,
\Big[ 1065 (B_{0, 12} - 2 B_{0, 14} + B_{0, 16}) \feq
$$ + 3872 R^2 (B_{0, 16} - 4 B_{0, 18} + 6 B_{0, 20} - 4 B_{0, 22} +
B_{0, 24}) \Big]^{1/2}~. $$
}
The above functions determine the control Cauchy problem
\rref{conca}. The numerical solution of this problem for many sample
values of $R$ yields a picture as in items (i)(ii), page \pageref{pagitem},
with a critical Reynolds number
\beq R_{*} \in (0.08, 0.09)~. \feq
In Boxes 1a-1d we consider the case $R=0.08$,
giving information on the following functions of time:
the quantity $(2 \pi)^{3/2} |u^1_k(t)|$ for the
wave vector \hbox{$k=(1,1,0)$};
the estimators $\Dd_3$ and $\ep_3$; the solution $\Rr_3$ of the
control Cauchy problem, which is global.
In Boxes 2a-2d we consider the analogous functions
in the case $R=0.09$, where $\Rr_3$ diverges at $\Tc = 2.153...$
({\footnote{Here and in the sequel, an
expression like $r= a. b c d e...~$ means that $a.b c d e$ are
the first digits of the MATHEMATICA output in the computation of $r$.}}).
Each one of these boxes (and of the subsequent ones) contains
the graph of the function under consideration, and its
numerical values for some choices of $t$. \par
One immediately notices that the functions in boxes of the types (a) and (b) (i.e., the norm
of the Fourier component $(1,1,0)$ and the estimator $\Dd_3$)
are very similar in these examples (and in all the
subsequent ones), even from the quantitative viewpoint. What really makes the difference
among these examples (and the forthcoming ones) are the
differential error estimator $\ep_3$ and the solution $\Rr_3$
of the control Cauchy problem, considered in type (c) and (d) boxes.
\salto
\textbf{Case $\boma{N=2}$.} This relies on the function
\beq u^2 := u_0 + R u_1 + R^2 u_2 : [0,+\infty) \vain \HM{5}~. \feq
Of course, $u_0, u_1$ are as in Eqs. \rref{u0} \rref{u1};
$u_2 = (u_{2, 1}, u_{2, 2}, u_{2, 3})$ has a more complicated
expression, and in order to save room we only report its first
component. This is
\beq u_{2, 1} =
{1 \over 9} B_{0, 2} \, (-e_{(-1, -1, 0)} - e_{(-1, 0, -1)} - e_{(1, 0, 1)} -
e_{(1, 1, 0)})  \feq
$$
+ {1 \over 3} B_{0, 4} \, (e_{(0, -2, 0)} - e_{(0, 0, -2)} - e_{(0, 0, 2)} + e_{(0, 2, 0)})  $$
$$ + {1 \over 12} B_{0, 6} \, (
-e_{(-3, -2, -1)} - e_{(-3, -1, -2)} + e_{(-2, -3, -1)} - e_{(-2, -1, -3)} -
  e_{(-1, -3, -2)} - e_{(-1, -2, -3)}  $$
$$ + 4 e_{(-1, -1, 0)} + 4 e_{(-1, 0, -1)} -
  8 e_{(0, -2, 0)} + 8 e_{(0, 0, -2)} + 8 e_{(0, 0, 2)} - 8 e_{(0, 2, 0)} $$
$$ +  4 e_{(1, 0, 1)} + 4 e_{(1, 1, 0)} - e_{(1, 2, 3)} - e_{(1, 3, 2)} -
  e_{(2, 1, 3)} + e_{(2, 3, 1)} - e_{(3, 1, 2)} - e_{(3, 2, 1)} )
$$
$$ + {1 \over 9} B_{0, 8} \,
(e_{(-3, -2, -1)} + e_{(-3, -1, -2)} - e_{(-2, -3, -1)} + e_{(-2, -1, -3)} +
  e_{(-1, -3, -2)} + e_{(-1, -2, -3)}  $$
$$  - 2 e_{(-1, -1, 0)} - 2 e_{(-1, 0, -1)} +
  3 e_{(0, -2, 0)} - 3 e_{(0, 0, -2)} - 3 e_{(0, 0, 2)} + 3 e_{(0, 2, 0)}  $$
  $$ - 2 e_{(1, 0, 1)} - 2 e_{(1, 1, 0)} + e_{(1, 2, 3)} + e_{(1, 3, 2)} +
  e_{(2, 1, 3)} - e_{(2, 3, 1)} + e_{(3, 1, 2)} + e_{(3, 2, 1)})  $$
$$ + {1 \over 36} B_{0, 14} \, (-e_{(-3, -2, -1)} - e_{(-3, -1, -2)} + e_{(-2, -3, -1)} -
        e_{(-2, -1, -3)} - e_{(-1, -3, -2)} - e_{(-1, -2, -3)} $$
$$ - e_{(1, 2, 3)} - e_{(1, 3, 2)} - e_{(2, 1, 3)} + e_{(2, 3, 1)} - e_{(3, 1, 2)} - e_{(3, 2, 1)})~. $$
The next step involves the differential error $e(u^2)$ (see again
\rref{eunn}) and the (time dependent) norms
\beq \Dd_3 := \| u^2 \|_3, \quad
\Dd_{4} := \| u^2 \|_{4}, \quad
\ep_3 := \| e(u^2) \|_3~. \label{normss} \feq
For example, one has
\parn
\vbox{
$$ \Dd_3 = {\sqrt{2} \over 3 \sqrt{3}} (2 \pi)^{3/2} \Big[\, 1296 \, B_{0, 4}
+ 288 R^2\,(-B_{0, 4} + 84 \, B_{0, 8} - 164 \, B_{0, 10} +
   81 B_{0, 12}) $$
$$ + R^4\,(16\,B_{0, 4} + 1056\,B_{0, 8} - 4544\,B_{0, 10}
+ 16317\,B_{0, 12} - 29496\,B_{0, 14} + 17680\,B_{0, 16} $$
\beq  + 6174\,B_{0, 20} -
   8232\,B_{0, 22} + 1029\,B_{0, 28}) \, \Big]^{1/2}~. \feq
}
The expressions of $\Dd_4$ and $\ep_3$ are not reported. The former
has a complexity similar to that of $\Dd_3$, the latter is
lengthier; in fact, $\ep_3$ has the form
$(2 \pi)^{3/2}
(\sum_{j, b} C_{j, b} R^j B_{0, b})^{1/2}$,
where $C_{j, b}$ are rational coefficients
and the sum involves $48$ pairs $(j,b)$, with $j \in \{6,8,10\}$.
\par
Let us pass to the control Cauchy problem
\rref{conca}; in the present case, the
picture of items (i)(ii), page \pageref{pagitem}
is realized with a critical Reynolds number
\beq R_{*} \in (0.12, 0.13)~. \label{1213} \feq
Boxes 3a-3d are about the case $R=0.12$,
and give information on the functions already
chosen for the previous tables; one of them is the solution $\Rr_3$ of the
control Cauchy problem, which is global.
Boxes 4a-4d are about the analogous functions
in the case $R=0.13$, in which $\Rr_3$ diverges at $\Tc = 2.604...$. \par
\salto
\textbf{Case $\boma{N=5}$.} This relies on the function
\beq u^5 := \sum_{j=0}^5 R^j u_j  : [0,+\infty) \vain \HM{5}~. \feq
The terms $u_0, u_1, u_2$ are as before;
the functions $u_3, u_4$ and $u_5$ have expressions of increasing complexity, that
we cannot reproduce here. We only mention that each one of the components
$u_{5, 1}$, $u_{5, 2}$, $u_{5, 3}$ is a linear combination with rational
coefficients of $1924$ terms of
the form $i B_{a, b} e_k$, with $a \in \{0,1,2\}$; the wave vectors $k$
appearing at least in one of the three components of $u_{5}$ are $174$.
\par
The differential error
and the norms
\beq \Dd_3 := \| u^5 \|_3, \quad
\Dd_{4} := \| u^5 \|_{4}, \quad
\ep_3 := \| e(u^5)\|_3 \label{norms} \feq
have very lengthy expressions.
Each one of the functions \rref{norms} has the form $(2 \pi)^{3/2}
(\sum_{j, a, b} C_{j, a, b} R^j B_{a, b})^{1/2}$,
where the $C_{j, a, b}$ are rational coefficients and
the sum involves finitely many
triples $(j, a, b)$ of nonnegative integers; these triples
are $204$ in the cases of $\Dd_3$ and $\Dd_4$, and $1734$ in the case
of $\ep_3$. \par
As in the other cases, the final step is the control problem \rref{conca};
the picture outlined by items items (i)(ii), page \pageref{pagitem}
is now realized with a critical Reynolds number
\beq R_{*} \in (0.23, 0.24)~. \feq
In Boxes 5a-5d and 6a-6d we give
some information about the cases $R=0.23$ and $R=0.24$, respectively.
Boxes 7a-7d are about the case $R = 0.12$; this choice
is considered as well in the next subsection, where it is used for a
comparison between the Reynolds expansion and the Galerkin approach.
\salto
\textbf{Comparison with the Galerkin approach.}
In \cite{appeul} we have considered the NS Cauchy problem
with the BNW initial datum, using the Galerkin approximate solution
$\ug$ that corresponds to a finite set of Fourier modes $G \subset \Zt \setminus \{0 \}$
(this is assumed to be reflection invariant: $k \in G$ $\Rightarrow$ $-k \in G$).
\par
This approach relies on the (finite-dimensional) Galerkin subspace
\beq \HG := \{ v \in \Dsz~|~v_k = 0 \quad \mbox{for}~k \not\in G \} \label{hg} \feq
$$ = \{ \sum_{k \in \G} v_k e_k~|~v_k \in \complessi^d,~\overline{v_k} = v_{-k}~,
k \sc \, v_k = 0~\mbox{for all $k \in G$} \} $$
and on the projection
\beq \EG : \Dsz~ \vain \HG~, \qquad v =
\sum_{k \in \Zd \setminus \{ 0 \}} v_k e_k \mapsto \EG v := \sum_{k \in \G} v_k e_k~. \label{pg} \feq
The Galerkin approximate solution of the Cauchy problem \rref{caun} (with the
BNW datum) corresponding
to the set of modes $G$ is the unique function $\ug$ such that
\beq \ug \in C^1([0,+\infty), \HG)~, \quad
{d \ug \over dt} = \Delta \ug + R \, \EG \PPP(\ug, \ug)~, \quad \ug(0) = \EG \uz
\label{galer} \feq
(thus, $\ug$ solves a finite-dimensional Cauchy problem; one can show the
existence of a global solution, of domain $[0,+\infty)$, using
the fact that the $L^2$ norm is a decreasing function of time).
Indeed, in \cite{appeul} we considered an equivalent formulation of
\rref{galer} based on the viscosity $\nu$ and on the unscaled time
$\ti$, related to $R$ and to the present time variable $t$ via Eq.
\rref{setting}; in the sequel we will rephrase the results of
\cite{appeul} in terms of the variables $R$ and $t$.
\par
For a given $G$, $\ug$ can be computed numerically; more precisely,
one solves numerically a system of ODEs for the Fourier coefficients
$\gaamma_k$ in the expansion
\beq \ug(t) = \sum_{k \in G} \gaamma_k(t) e_k~. \feq
One can specialize to $\ug$ the general framework
for approximate NS solutions; in particular, one introduces
the differential error $e(\ug) := d \ug/d t - \Delta \ug - R \PPP(\ug, \ug)$
and the tautological estimators of orders 3 or 4 for the growth of $\ug$ and for the above
error, i.e., the functions
\beq \Ddd_3(t) := \| \ug(t) \|_3,~\Ddd_4(t) := \| \ug(t) \|_4,~
\eep_3(t) := \| e(\ug)(t) \|_3
\quad (t \in [0,+\infty)) \feq
(see \cite{appeul} for an explicit expression of $\eep_3$ in terms
of the Fourier coefficients $\gaamma_k$). The datum error
$\ug(0) - \uz$ is zero if $G$ contains the Fourier
modes $\pm a, \pm b, \pm c$ involved in Eq. \rref{unec}.
Assuming this, one can analyze the Galerkin
approximate solution in terms of a control Cauchy problem
\beq {d \Rrr_3 \over d t} = - \Rrr_3
+ R (G_3 \Ddd_3 + K_3 \Ddd_{4}) \Rrr_3 + R \, G_3 \Rrr^2_3 + \eep_3
~\mbox{on $[0,\Tcc)$},~~\Rrr_3(0) =0~, \label{concauu} \feq
whose (maximal) solution $\Rrr_n \in C^1([0,\Tcc), [0,+\infty))$
can be computed numerically.
In \cite{appeul} we have employed a set $G$ of $150$ modes, of the following form:
\par
\vbox{
\beq G := S \cup -S~; \qquad -S := \{ -k~|~k \in S \}~; \label{defg} \feq
$$ S :=
\{ (0, 0, 2), (0, 1, -3), (0, 1, 1), (0, 1, 3), (0, 2, 0), (0, 2, 2), (0, 3, -1),
(0, 3, 1), (0, 3, 3), $$
$$ (1, -3, -2), (1, -3, 0), (1, -3, 2), (1, -2, -3),
(1, -2, -1), (1, -2, 1), (1, -2, 3), $$
$$  (1, -1, -2), (1, -1, 2), (1, 0, -3),
(1, 0, 1), (1, 0, 3),
(1, 1, -2), (1, 1, 0),
(1, 1, 2), (1, 2, -3), $$
$$ (1, 2, -1), (1, 2, 1), (1, 2, 3), (1, 3, -2), (1, 3, 0),
(1, 3, 2), (2, -3, -3),
(2, -3, -1), (2, -3, 1), $$
$$ (2, -3, 3), (2, -2, -2), (2, -2, 2), (2, -1, -3),
(2, -1, -1), (2, -1, 1),
(2, -1, 3), (2, 0, 0), $$
$$ (2, 0, 2), (2, 1, -3), (2, 1, -1), (2, 1, 1), (2, 1, 3),
(2, 2, -2), (2, 2, 0),
(2, 3, -3), (2, 3, -1), $$
$$ (2, 3, 1), (2, 3, 3), (3, -3, -2), (3, -3, 2), (3, -2, -3),
(3, -2, -1),
(3, -2, 1), (3, -2, 3), $$
$$ (3, -1, -2), (3, -1, 0), (3, -1, 2), (3, 0, -1), (3, 0, 1),
(3, 0, 3), (3, 1, -2),
(3, 1, 0), (3, 1, 2), $$
$$ (3, 2, -3), (3, 2, -1), (3, 2, 1), (3, 2, 3), (3, 3, -2),
(3, 3, 0), (3, 3, 2)  \}~.
$$
}
\noindent
Let us summarize the results arising from this choice of $G$.
\begin{itemize}
\item[(\,$\hat{\imath}$\,)] There are indications for
the existence of a critical Reynolds number $\Rc$ such that:
$\Tcc = +\infty$ and $\Rrr_3$ vanishes at $+\infty$
if $0 \leqs R \leqs \Rc$, while $\Tcc < +\infty$ and $\Rrr_3$ diverges at
$\Tcc$ if $R > \Rc$. One has $\Rc \in (0.13, 0.14)$
({\footnote{The estimate coming directly from \cite{appeul}
is $1/8 < \Rc < 1/7$, whence $0.125 < \Rc < 0.143$;
this estimate has been refined to $0.13 < \Rc < 0.14$ by
a supplementary run of the MATHEMATICA program for \cite{appeul}.}}).
\item[(\,$\widehat{\imath \imath}$\,)] Let us consider the maximal solution $u$ of the NS Cauchy problem \rref{caun}.
For $0 \leqs R \leqs \Rc$, (\,$\hat{\imath}$\,) grants that $u$ is global and $\| u(t) - \ug(t) \|_3
\leqs \Rrr_3(t)$ for all $t \in [0,+\infty)$. For $R > \Rc$
(\,$\hat{\imath}$\,) only grants that the domain of $u$ contains the interval $[0, \Tcc)$, and
that $\| u(t) - \ug(t) \|_3 \leqs \Rrr_3(t)$ for all $t \in [0,\Tcc)$.
\end{itemize}
As an example, Boxes 8a-8d report the main results about the Galerkin approach for
$R=0.12$.
Both the general picture ($\widehat{\imath}$) ($\widehat{\imath \imath}$) and the results for $R=0.12$
indicate a substantial equivalence between the Galerkin approach with the set
$G$ of \rref{defg} and the Reynolds expansion of order $N=2$.
In fact:
\begin{itemize}
\item[(a)] The expansion of order $N=2$ yields a picture
similar to ($\widehat{\imath}$) ($\widehat{\imath \imath}$), with
a critical Reynolds number $R_{*} \in (0.12, 0.13)$ (see
items (i)(ii) of page \pageref{pagitem} and Eq. \rref{1213} for $R_{*}$).
\item[(b)] The error estimators $\ep_3$, $\eep_3$
and the solutions $\Rr_3$, $\Rrr_3$ of
the control equations for the $N=2$ Reynolds expansion
and for the Galerkin approach with $G$ in \rref{defg}
are have similar orders of magnitude in the case
$R=0.12$, here used as a test to make comparisons:
see Boxes 3c-3d and 8c-8d for more detailed information.
\end{itemize}
If we pass to the $N=5$ expansion, we find
a significant improvement with respect to the results of
the above Galerkin approach. Let us recall that, for $N=5$,
the critical Reynolds number for global existence
is in the interval $(0.23, 0.24)$; moreover, if we use again
the case $R=0.12$ for a comparison, we see that
the error estimators $\ep_3$ and the  function $\Rr_3$
of the $N=5$ expansion are much smaller than
the homologous function $\eep_3, \Rrr_3$ of the Galerkin approach: see
Boxes 7c-7d and 8c-8d. In particular, the ratio
$\Rr_3(t)/\Rrr_3(t)$ is of order $10^{-4}$:
so, the bound $\| u(t) - u^5(t) \|_3 \leqs \Rr_3(t)$ is much
more stringent than the bound $\| u(t) - \ug(t) \|_3 \leqs
\Rrr_3(t)$.
\begin{figure}
\framebox{
\parbox{2in}{
\includegraphics[
height=1.3in,
width=2.0in
]%
{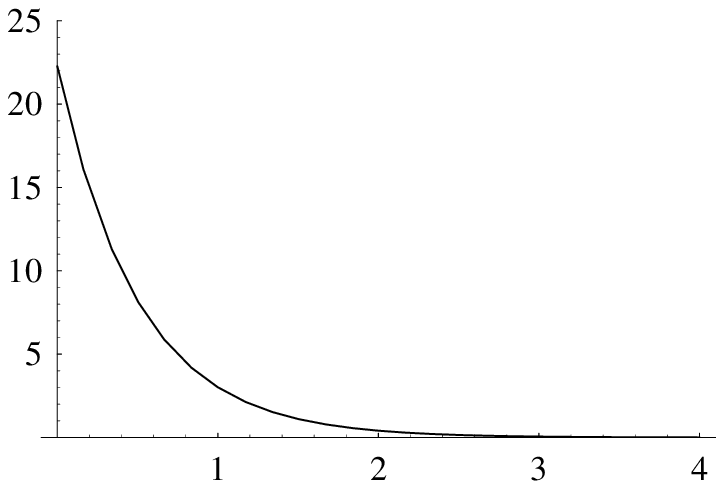}%
\par
{\tiny{
{\textbf{Box 1a.}
$N=1$, $R = 0.08$:
the function $\gamma(t) := (2 \pi)^{3/2} |u^1_{(1,1,0)}(t)|$. One has
$\gamma(0) = 22.27...$, $\gamma(0.5) = 8.189...$,
$\gamma(1) = 3.012...$, $\gamma(1.5) = 1.108... $, $\gamma(2) = 0.4076...$,
$\gamma(4) = 7.466... \times 10^{-3}$,
$\gamma(8) = 2.504... \times 10^{-6}$, $\gamma(10) = 4.587... \times 10^{-8}$~.
}}
\par}
\label{f1a}
}
}
\hskip 0.4cm
\framebox{
\parbox{2in}{
\includegraphics[
height=1.3in,
width=2.0in
]%
{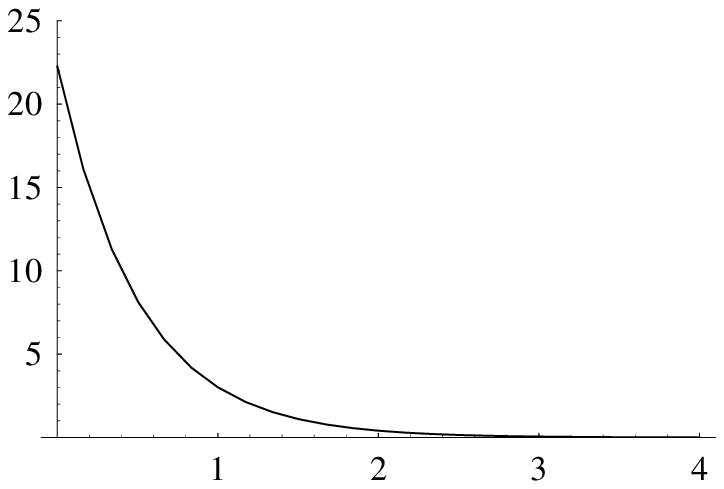}%
\par
{\tiny{
{\textbf{Box 2a.~}
$N=1$, $R=0.09$: the function
$\gamma(t) := (2 \pi)^{3/2} |u^1_{(1,1,0)}(t)|$. One has
$\gamma(0) = 22.27...$, $\gamma(0.5) = 8.188...$,
$\gamma(1) = 3.011...$, $\gamma(1.5) = 1.107... $, $\gamma(2) = 0.4075...$,
$\gamma(4) = 7.465... \times 10^{-3}$~.
{~}
\vskip 0.05cm \noindent
}}
\par}
\label{f2a}
}
}
\framebox{
\parbox{2in}{
\includegraphics[
height=1.3in,
width=2.0in
]%
{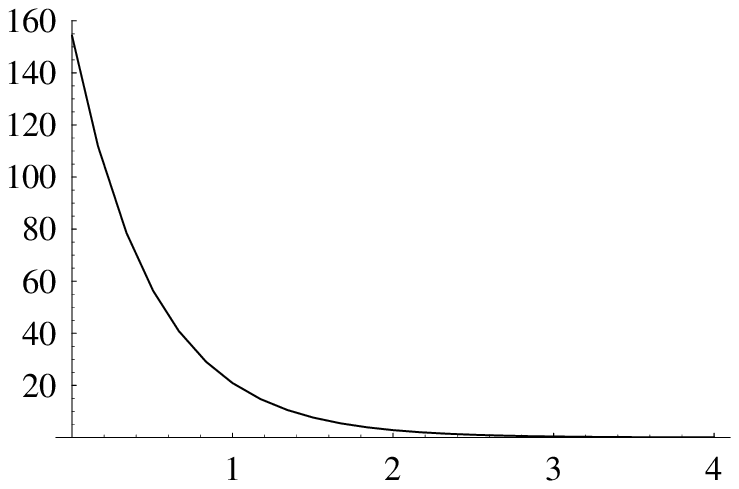}%
\par
{\tiny{
{\textbf{Box 1b.}
$N=1, R=0.08$: the \hbox{function}
\hbox{$\Dd_3(t) := \| u^1(t) \|_3$.
One has $\Dd_3(0) = 154.3...$}, $\Dd_3(0.5) = 56.94...$,
$\Dd_3(1) = 20.90...$, $\Dd_3(1.5) = 7.683... $, $\Dd_3(2) = 2.826...$,
$\Dd_3(4) = 0.05176...$, $\Dd_3(8) =1.736... \times 10^{-5}$,
$\Dd_3(10) = 3.180... \times 10^{-7}$~.
}}
\par}
\label{f1b}
}
}
\hskip 0.4cm
\framebox{
\parbox{2in}{
\includegraphics[
height=1.3in,
width=2.0in
]%
{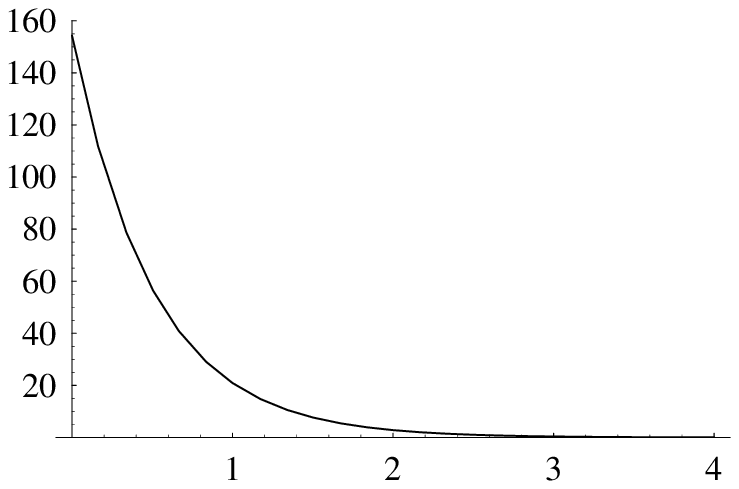}%
\par
{\tiny{
{\textbf{Box 2b.~}
$N=1$, $R=0.09$: \hbox{the function}
\hbox{$\Dd_3(t) := \| u^1(t) \|_3$.
One has $\Dd_3(0) = 154.3...$}, $\Dd_3(0.5) = 56.99...$,
$\Dd_3(1) = 20.90...$, $\Dd_3(1.5) = 7.684... $, $\Dd_3(2) = 2.826...$,
$\Dd_3(4) = 0.05176...$~.
{~}
\vskip 0.2cm \noindent
}}
\par}
\label{f2b}
}
}
\framebox{
\parbox{2in}{
\includegraphics[
height=1.3in,
width=2.0in
]%
{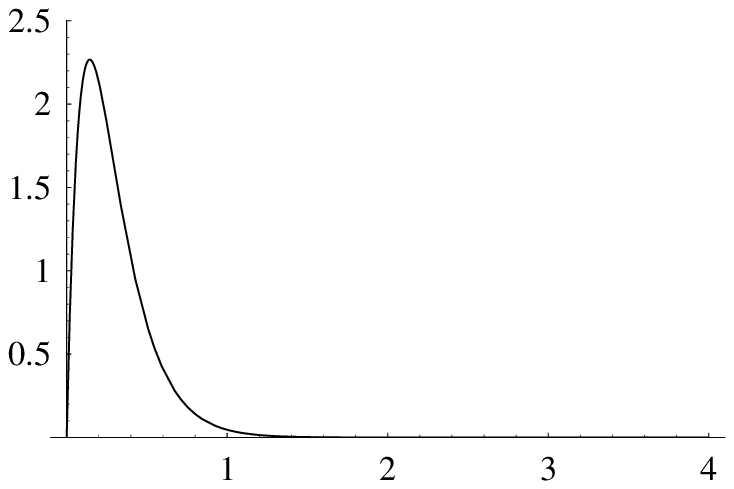}%
\par
{\tiny{
{\textbf{Box 1c.~}
$N=1$, $R=0.08$: the function
$\ep_3(t)$.
One has $\ep_3(0) = 0$, $\ep_3(0.14) = 2.266...$,
$\ep_3(1) = 0.04605...$,
$\ep_3(2) = 1.296... \times 10^{-4}$, $\ep_3(4) = 8.108... \times 10^{-10}$\, .
}}
\par}
\label{f1c}
}
}
\hskip 0.4cm
\framebox{
\parbox{2in}{
\includegraphics[
height=1.3in,
width=2.0in
]%
{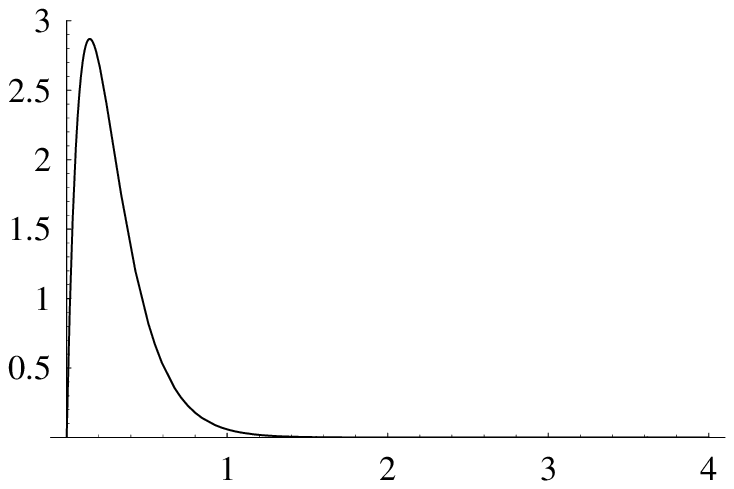}%
\par
{\tiny{
{\textbf{Box 2c.~}
$N=1$, $R=0.09$: the function $\ep_3(t)$.
One has $\ep_3(0) = 0$,
$\ep_3(0.14) = 2.868...$,
$\ep_3(1) = 0.05829...$,
$\ep_3(2) = 1.640... \times 10^{-4}$,
$\ep_3(4) = 1.026... \times 10^{-9}$\, .
\vskip 0.05cm
{~}
}}
\par}
\label{f2c}
}
}
\framebox{
\parbox{2in}{
\includegraphics[
height=1.3in,
width=2.0in
]%
{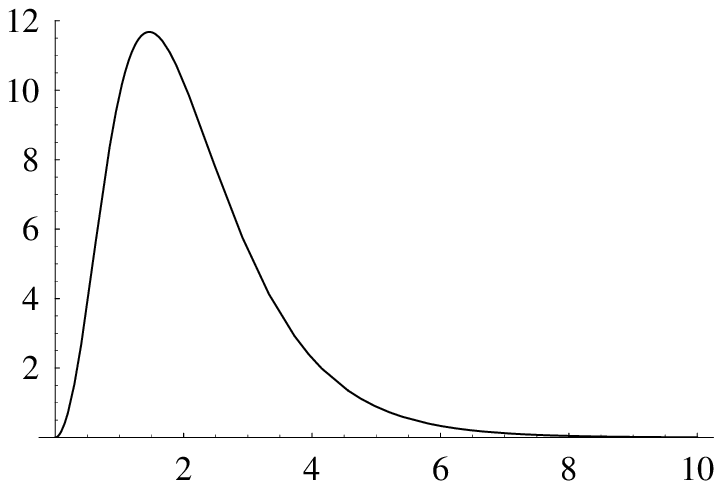}%
\par
{\tiny{
{\textbf{Box 1d.~}
$N=1$, $R=0.08$: the function
$\Rr_3(t)$. This
appears to be \hbox{globally} defined, and vanishing at $+\infty$.
One has $\Rr_3(0) = 0$, $\Rr_3(1) = 9.858...$,
$\Rr_3(1.5) = 11.66...$, $\Rr_3(2) = 10.23...$, $\Rr_3(4) = 2.283...$,
$\Rr_3(8) = 0.04547...$, $\Rr_3(10) = 6.162... \times 10^{-3}$~.
\vskip 0.1cm
{~}
}}
\par}
\label{f1d}
}
}
\hskip 3.7cm
\framebox{
\parbox{2in}{
\includegraphics[
height=1.3in,
width=2.0in
]%
{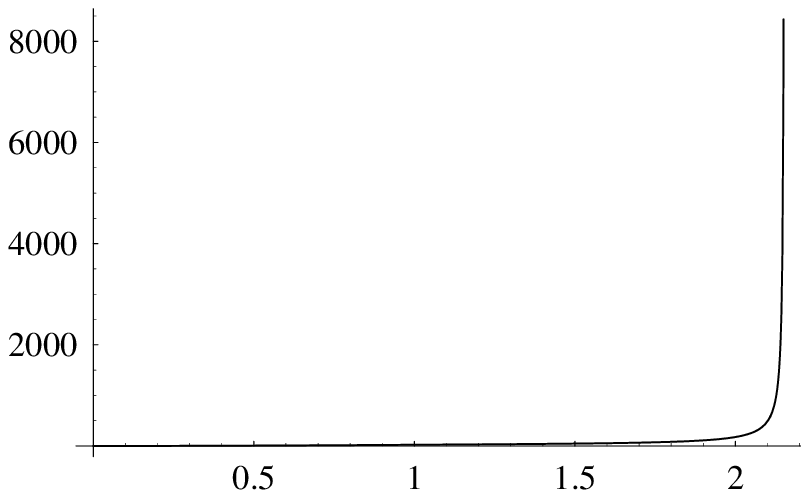}%
\par
{\tiny{
{\textbf{Box 2d.~}
\hbox{$N=1$, $R=0.09$:
the function} $\Rr_3(t)$. This diverges
for $t \vain \Tc = 2.153...$\,.
One has $\Rr_3(0) = 0$, $\Rr_3(0.5) = 6.624...$,
$\Rr_3(1) = 23.20...$, $\Rr_3(1.5) = 46.60...$,
$\Rr_3(2) = 176.0...$~.
\vskip 0.35cm
{~}
}}
\par}
\label{f2d}
}
}
\end{figure}
\begin{figure}
\framebox{
\parbox{2in}{
\includegraphics[
height=1.3in,
width=2.0in
]%
{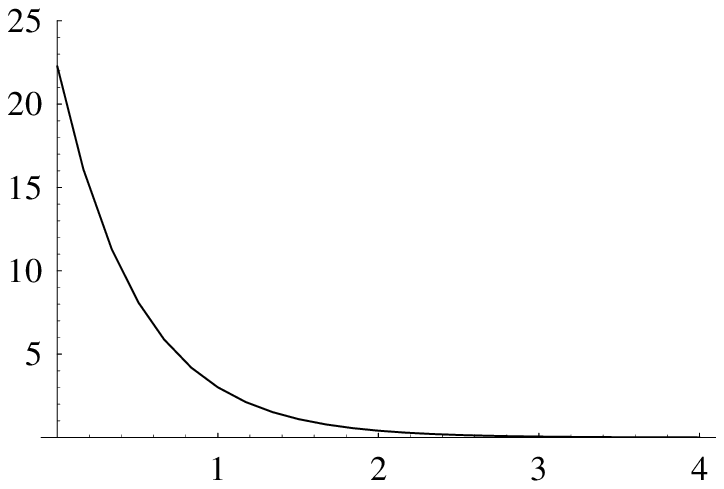}%
\par
{\tiny{
{\textbf{Box 3a.~}
$N=2$, $R = 0.12$:
the function $\gamma(t) := (2 \pi)^{3/2} |u^1_{(1,1,0)}(t)|$. One has
$\gamma(0) = 22.27...$, $\gamma(0.5) = 8.184...$,
$\gamma(1) = 3.009...$, $\gamma(1.5) = 1.107... $, $\gamma(2) = 0.4072...$,
$\gamma(4) = 7.459... \times 10^{-3}$,
$\gamma(8) = 2.502... \times 10^{-6}$, $\gamma(10) = 4.583... \times 10^{-8}$~.
}}
\par}
\label{f3a}
}
}
\hskip 0.4cm
\framebox{
\parbox{2in}{
\includegraphics[
height=1.3in,
width=2.0in
]%
{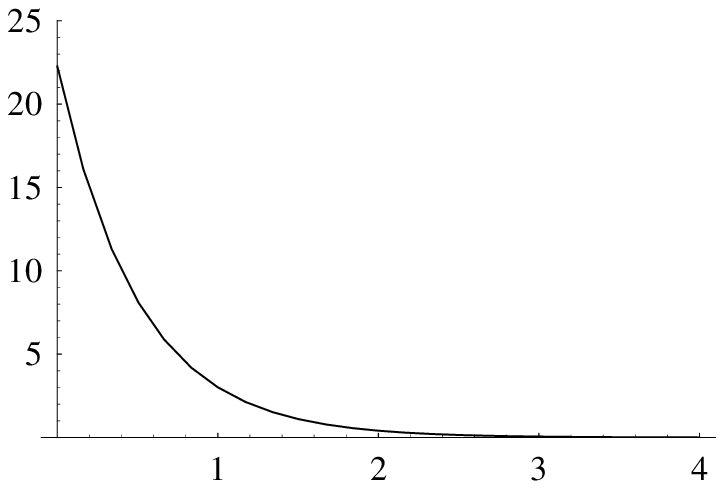}%
\par
{\tiny{
{\textbf{Box 4a.~}
$N=2$, $R=0.13$: the function
$\gamma(t) := (2 \pi)^{3/2} |u^1_{(1,1,0)}(t)|$. One has
$\gamma(0) = 22.27...$, $\gamma(0.5) = 8.183...$,
$\gamma(1) = 3.009...$, $\gamma(1.5) = 1.106... $, $\gamma(2) = 0.4071...$,
$\gamma(4) = 7.457... \times 10^{-3}$~.
\vskip 0.2cm
{~}
}}
\par}
\label{f4a}
}
}
\framebox{
\parbox{2in}{
\includegraphics[
height=1.3in,
width=2.0in
]%
{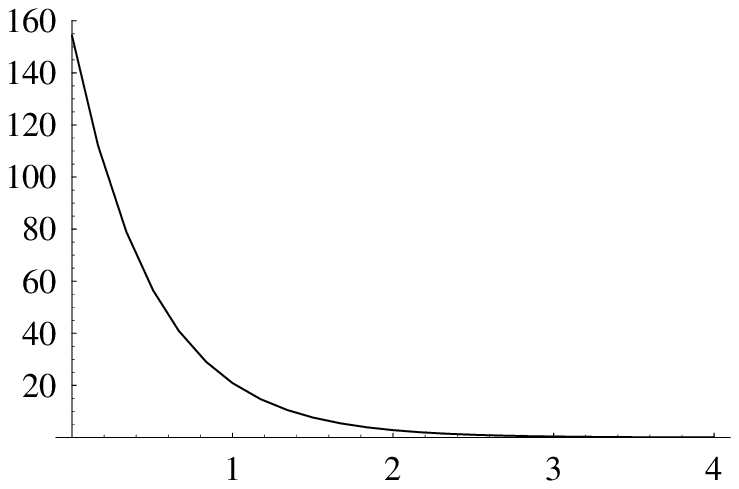}%
\par
{\tiny{
{\textbf{Box 3b.~}
\hbox{$N=2$, $R=0.12$: the function}
\hbox{$\Dd_3(t) := \| u^1(t) \|_3$.
One has $\Dd_3(0) = 154.3...$}, $\Dd_3(0.5) = 57.10...$,
$\Dd_3(1) = 20.88...$, $\Dd_3(1.5) = 7.672... $, $\Dd_3(2) = 2.821...$,
$\Dd_3(4) = 0.05168...$, $\Dd_3(8) =1.733... \times 10^{-5}$,
$\Dd_3(10) = 3.175... \times 10^{-7}$~.
}}
\par}
\label{f3b}
}
}
\hskip 0.4cm
\framebox{
\parbox{2in}{
\includegraphics[
height=1.3in,
width=2.0in
]%
{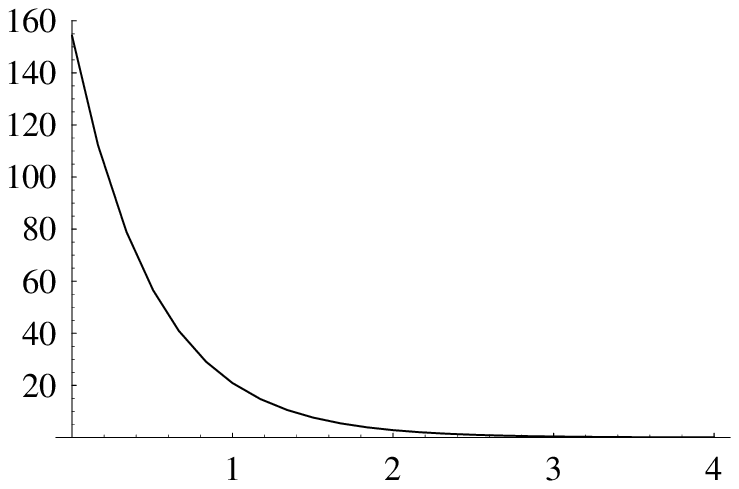}%
\par
{\tiny{
{\textbf{Box 4b.~}
\hbox{$N=2$, $R=0.13$: the function}
\hbox{$\Dd_3(t) := \| u^1(t) \|_3$.
One has $\Dd_3(0) = 154.3...$}, $\Dd_3(0.5) = 57.16...$,
$\Dd_3(1) = 20.89...$, $\Dd_3(1.5) = 7.671... $, $\Dd_3(2) = 2.821...$,
$\Dd_3(4) = 0.05166...$~.
\vskip 0.15cm
{~}
}}
\par}
\label{f4b}
}
}
\framebox{
\parbox{2in}{
\includegraphics[
height=1.3in,
width=2.0in
]%
{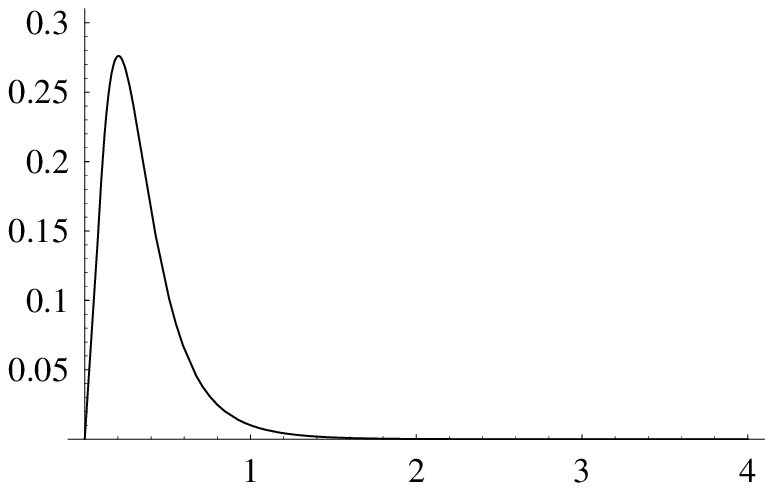}%
\par
{\tiny{
{\textbf{Box 3c.~}
\hbox{$N=2$, $R=0.12$: the function}
\hbox{$\ep_3(t)$.
One has $\ep_3(0) = 0$, $\ep_3(0.20) = 0.2759...$,}
$\ep_3(1) = 0.01007...$,
$\ep_3(2) = 1.688... \times 10^{-4}$, $\ep_3(4) = 5.653... \times 10^{-8}$\, .
}}
\par}
\label{f3c}
}
}
\hskip 0.4cm
\framebox{
\parbox{2in}{
\includegraphics[
height=1.3in,
width=2.0in
]%
{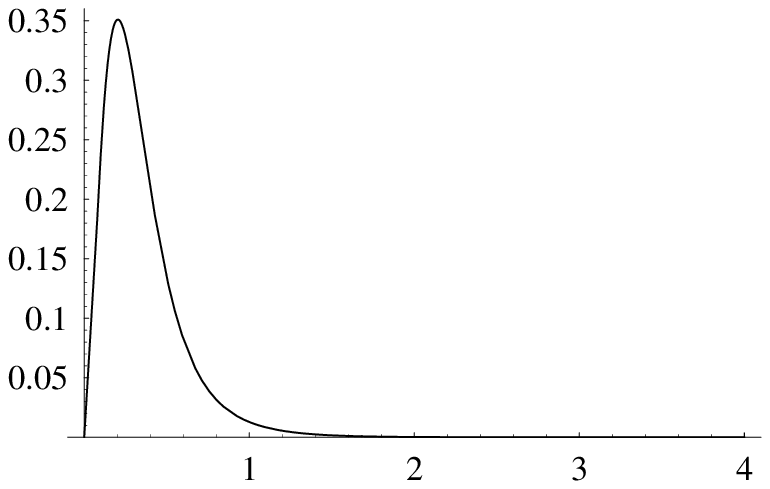}%
\par
{\tiny{
{\textbf{Box 4c.~}
\hbox{$N=2$, $R=0.13$: the function}
\hbox{$\ep_3(t)$. One has $\ep_3(0) = 0$,
$\ep_3(0.20) = 0.3508...$},
$\ep_3(1) = 0.01280...$,
$\ep_3(2) = 2.146... \times 10^{-4}$,
$\ep_3(4) = 7.187... \times 10^{-8}$\, .}}
\par}
\label{f4c}
}
}
\framebox{
\parbox{2in}{
\includegraphics[
height=1.3in,
width=2.0in
]%
{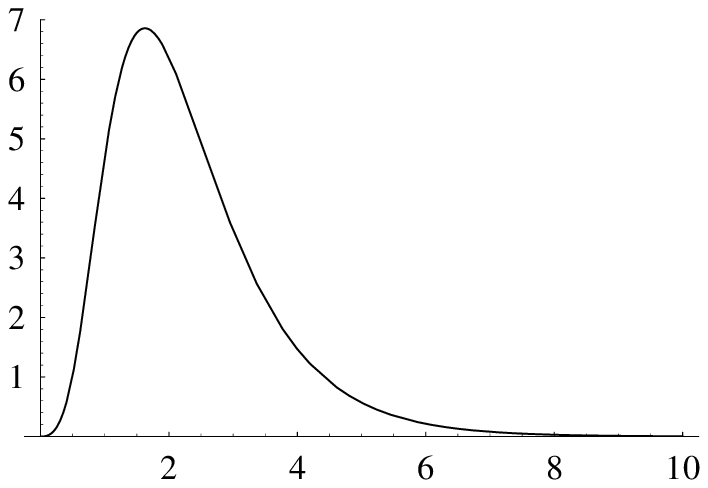}%
\par
{\tiny{
{\textbf{Box 3d.~}
\hbox{$N=2$, $R=0.12$: the function}
\hbox{$\Rr_3(t)$. This
appears to be globally defined,}
\hbox{and vanishing at $+\infty$. One has $\Rr_3(0) = 0$}, $\Rr_3(1) = 4.699...$,
$\Rr_3(1.6) = 6.854...$, $\Rr_3(2) = 6.368...$, $\Rr_3(4) = 1.466...$,
$\Rr_3(8) = 0.02914...$, $\Rr_3(10) = 3.949... \times 10^{-3}$~.
\vskip 0.1cm
{~}
}}
\par}
\label{f3d}
}
}
\hskip 3.7cm
\framebox{
\parbox{2in}{
\includegraphics[
height=1.3in,
width=2.0in
]%
{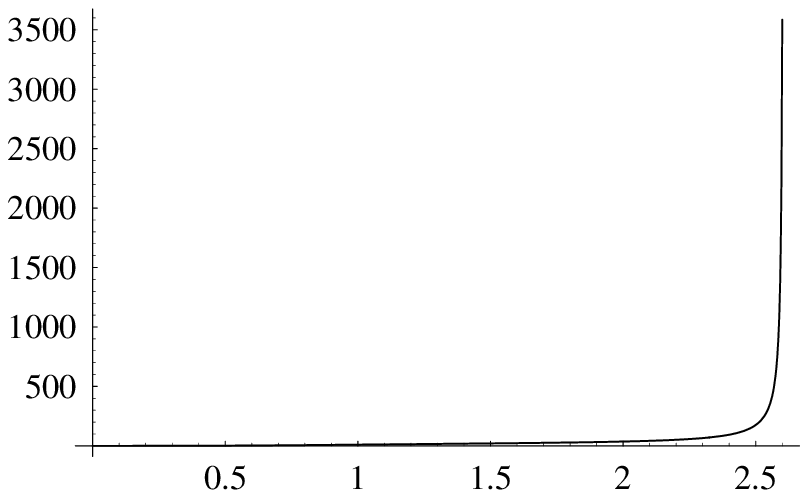}%
\par
{\tiny{
{\textbf{Box 4d.~}
\hbox{$N=2$, $R=0.13$: the function} $\Rr_3(t)$.
\hbox{This diverges
for $t \vain \Tc = 2.604...$\,.}
One has $\Rr_3(0) = 0$, $\Rr_3(0.5) = 1.735...$,
$\Rr_3(1) = 10.20...$, $\Rr_3(1.5) = 20.68...$,
$\Rr_3(2) = 36.30...$, $\Rr_3(2.5) = 175.3$~.
\vskip 0.25cm
{~}
}}
\par}
\label{f4d}
}
}
\end{figure}
\begin{figure}
\framebox{
\parbox{2in}{
\includegraphics[
height=1.3in,
width=2.0in
]%
{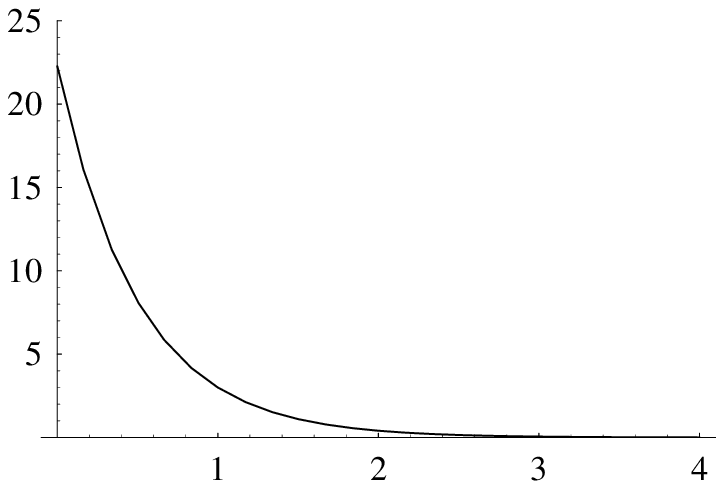}%
\par
{\tiny{
{\textbf{Box 5a.~}
\hbox{$N=5$, $R = 0.23$: the function}
$\gamma(t) := (2 \pi)^{3/2} |u^5_{(1,1,0)}(t)|$. One has
$\gamma(0) = 22.27...$, $\gamma(0.5) = 8.160...$,
$\gamma(1) = 2.997...$, $\gamma(1.5) = 1.102... $, $\gamma(2) = 0.4055...$,
$\gamma(4) = 7.428... \times 10^{-3}$,
$\gamma(8) = 2.491... \times 10^{-6}$, $\gamma(10) = 4.564... \times 10^{-8}$~.
}}
\par}
\label{f5a}
}
}
\hskip 0.4cm
\framebox{
\parbox{2in}{
\includegraphics[
height=1.3in,
width=2.0in
]%
{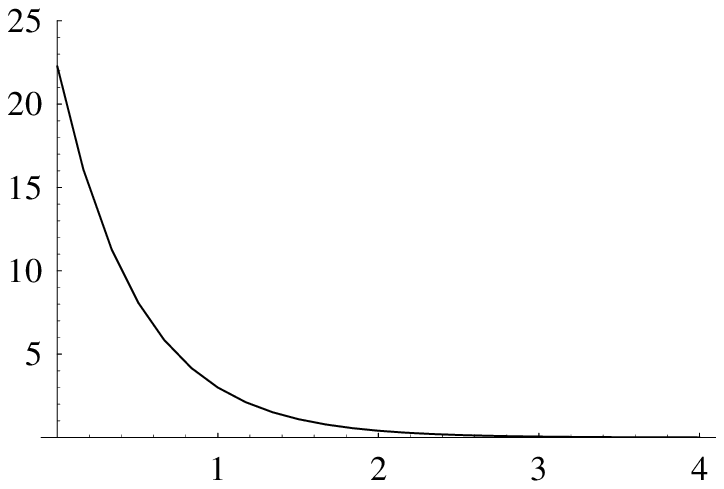}%
\par
{\tiny{
{\textbf{Box 6a.~}
\hbox{$N=5$, $R=0.24$: the function}
$\gamma(t) := (2 \pi)^{3/2} |u^5_{(1,1,0)}(t)|$. One has
$\gamma(0) = 22.27...$, $\gamma(0.5) = 8.157...$,
$\gamma(1) = 2.996...$, $\gamma(1.5) = 1.101... $, $\gamma(2) = 0.4053...$,
$\gamma(3) = 0.05485...$,
$\gamma(4) = 7.424... \times 10^{-3}$~.
\vskip 0.15cm
{~}
}}
\par}
\label{f6a}
}
}
\framebox{
\parbox{2in}{
\includegraphics[
height=1.3in,
width=2.0in
]%
{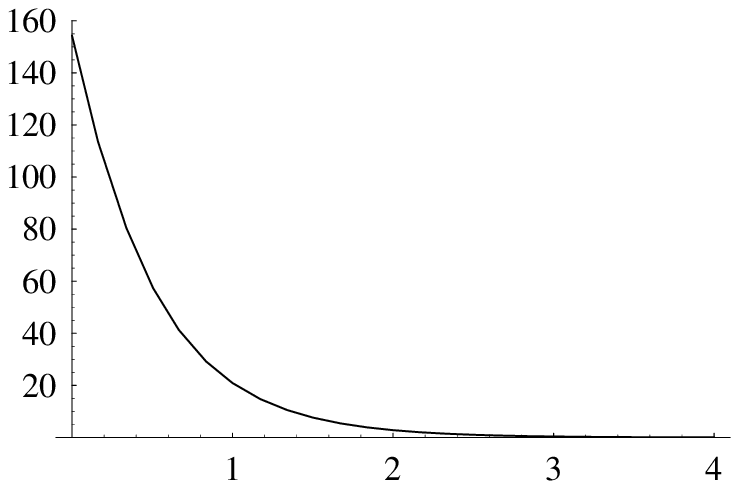}%
\par
{\tiny{
{\textbf{Box 5b.~}
$N=5$, $R=0.23$: the function
$\Dd_3(t) = \| u^5(t) \|_3$.
One has $\Dd_3(0) = 154.3...$, $\Dd_3(0.5) = 56.97...$,
$\Dd_3(1) = 20.90...$, $\Dd_3(1.5) = 7.646... $, $\Dd_3(2) = 2.810...$,
$\Dd_3(4) = 0.05146...$, $\Dd_3(8) =1.726... \times 10^{-5}$,
$\Dd_3(10) = 3.162... \times 10^{-7}$~.
}}
\par}
\label{f5b}
}
}
\hskip 0.4cm
\framebox{
\parbox{2in}{
\includegraphics[
height=1.3in,
width=2.0in
]%
{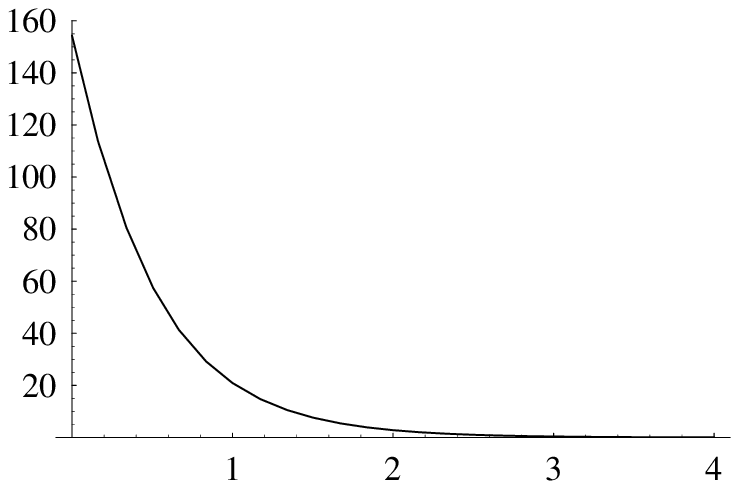}%
\par
{\tiny{
{\textbf{Box 6b.~}
$N=5$, $R=0.24$: the function
$\Dd_3(t) = \| u^5(t) \|_3$.
One has $\Dd_3(0) = 154.3...$, $\Dd_3(0.5) = 58.08...$,
$\Dd_3(1) = 20.90...$, $\Dd_3(1.5) = 7.643... $, $\Dd_3(2) = 2.808...$,
$\Dd_3(3) = 0.3800$,
$\Dd_3(4) = 0.05143...$~.
\vskip 0.15cm
{~}
}}
\par}
\label{f6b}
}
}
\framebox{
\parbox{2in}{
\includegraphics[
height=1.3in,
width=2.0in
]%
{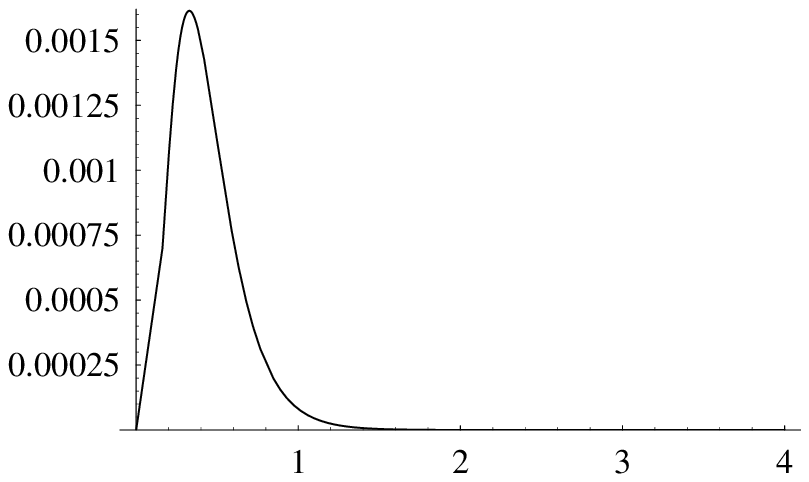}%
\par
{\tiny{
{\textbf{Box 5c.~}
\hbox{$N=5$, $R=0.23$: the function}
$\ep_3(t)$. One has $\ep_3(0) = 0$,
$\ep_3(0.33) = 1.613... \times 10^{-3}$,
$\ep_3(2) = 2.061... \times 10^{-7}$,
$\ep_3(3) = 8.234 \times 10^{-10}$,
$\ep_3(4) = 1.182... \times 10^{-11}$\, .
}}
\par}
\label{f5c}
}
}
\hskip 0.4cm
\framebox{
\parbox{2in}{
\includegraphics[
height=1.3in,
width=2.0in
]%
{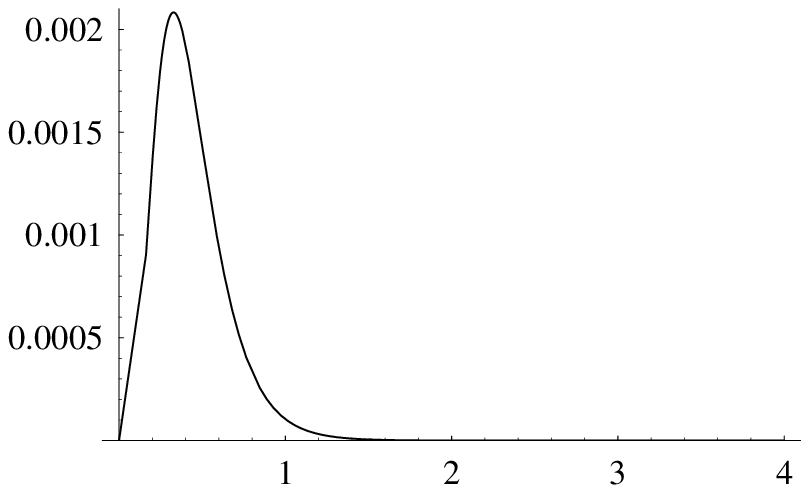}%
\par
{\tiny{
{\textbf{Box 6c.~}
\hbox{$N=5$, $R=0.24$: the function}
$\ep_3(t)$. One has $\ep_3(0) = 0$,
$\ep_3(0.33) = 2.0827 \times 10^{-3}...$,
$\ep_3(1) = 1.043 \times 10^{-4}...$,
$\ep_3(2) = 2.664... \times 10^{-7}$,
$\ep_3(4) = 1.591... \times 10^{-11}$\, .
\vskip 0.05cm
{~}
}}
\par}
\label{f6c}
}
}
\framebox{
\parbox{2in}{
\includegraphics[
height=1.3in,
width=2.0in
]%
{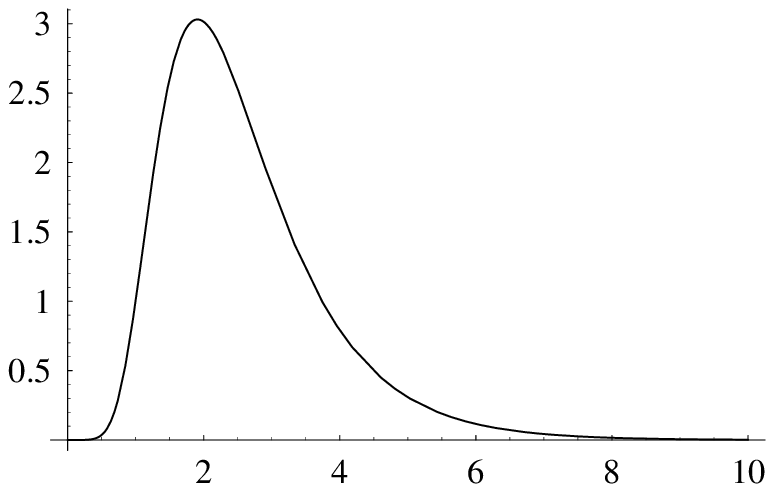}%
\par
{\tiny{
{\textbf{Box 5d.~}
$N=5$, $R=0.23$: the function $\Rr_3(t)$. This
appears to be \hbox{globally} defined, and vanishing at $+\infty$.
One has $\Rr_3(0) = 0$, $\Rr_3(1) = 1.004...$,
$\Rr_3(1.5) = 2.609...$, $\Rr_3(2) = 3.014...$, $\Rr_3(4) = 0.7907...$,
$\Rr_3(8) = 0.01580...$, $\Rr_3(10) = 2.141... \times 10^{-3}$~.
}}
\par}
\label{f5d}
}
}
\hskip 3.7cm
\framebox{
\parbox{2in}{
\includegraphics[
height=1.3in,
width=2.0in
]%
{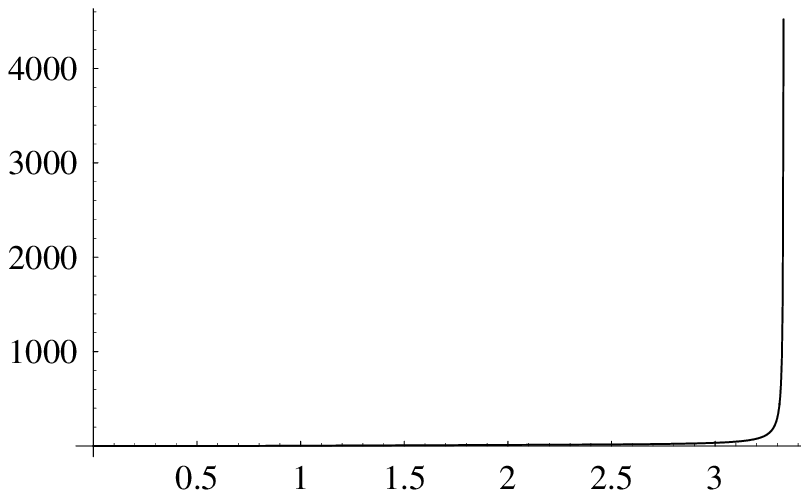}%
\par
{\tiny{
{\textbf{Box 6d.~}
$N=5$, $R=0.24$: the function $\Rr_3(t)$. This diverges
for $t \vain \Tc = 3.332...$\,.
One has $\Rr_3(0) = 0$, $\Rr_3(0.5) = 0.06348...$,
$\Rr_3(1) = 2.126...$, $\Rr_3(2) = 10.89...$,
$\Rr_3(3) = 33.31...$, $\Rr_3(3.33) = 4520.7...$~.
\vskip 0.25cm
{~}
}}
\par}
\label{f6d}
}
}
\end{figure}
\begin{figure}
\framebox{
\parbox{2in}{
\includegraphics[
height=1.3in,
width=2.0in
]%
{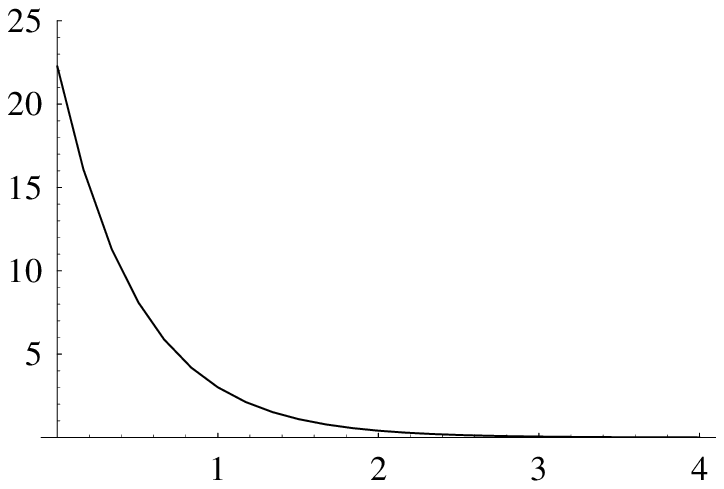}%
\par
{\tiny{
{\textbf{Box 7a.~}
$N=5$, $R = 0.12$: the function
$\gamma(t) := (2 \pi)^{3/2} |u^1_{(1,1,0)}(t)|$.
One has
$\gamma(0) = 22.27...$, $\gamma(0.5) = 8.184...$,
$\gamma(1) = 3.009...$, $\gamma(1.5) = 1.107... $, $\gamma(2) = 0.4072...$,
$\gamma(4) = 7.459... \times 10^{-3}$,
$\gamma(8) = 2.502... \times 10^{-6}$, $\gamma(10) = 4.583... \times 10^{-8}$~.
}}
\par}
\label{f7a}
}
}
\hskip 0.4cm
\framebox{
\parbox{2in}{
\includegraphics[
height=1.3in,
width=2.0in
]%
{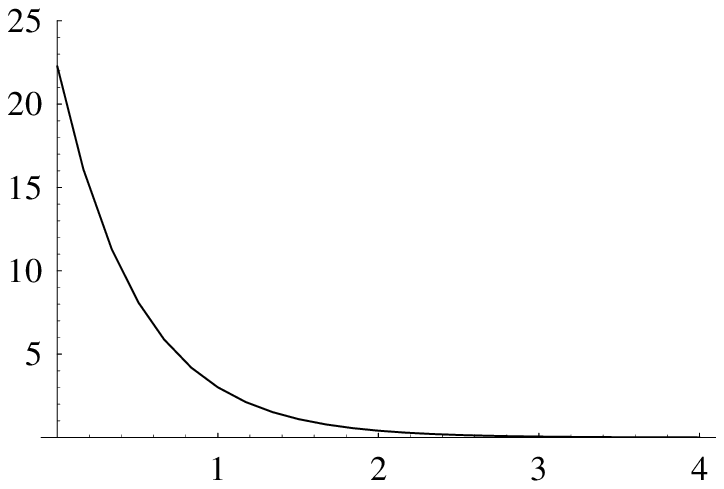}%
\par
{\tiny{
{\textbf{Box 8a.~}
Galerkin method with $G$ as in
\rref{defg}, for $R=0.12$: the function
$\gaamma(t) := (2 \pi)^{3/2} |\ug_{(1,1,0)}(t)|$.
One has
$\gaamma(0) = 22.27...$, $\gaamma(0.5) = 8.184...$,
$\gaamma(1) = 3.009...$, $\gaamma(1.5) = 1.107... $, $\gaamma(2) = 0.4072...$,
$\gaamma(4) = 7.459... \times 10^{-3}$,
$\gaamma(8) = 2.506... \times 10^{-6}$, $\gaamma(10) = 4.562... \times 10^{-8}$~.
}}
\par}
\label{f8a}
}
}
\framebox{
\parbox{2in}{
\includegraphics[
height=1.3in,
width=2.0in
]%
{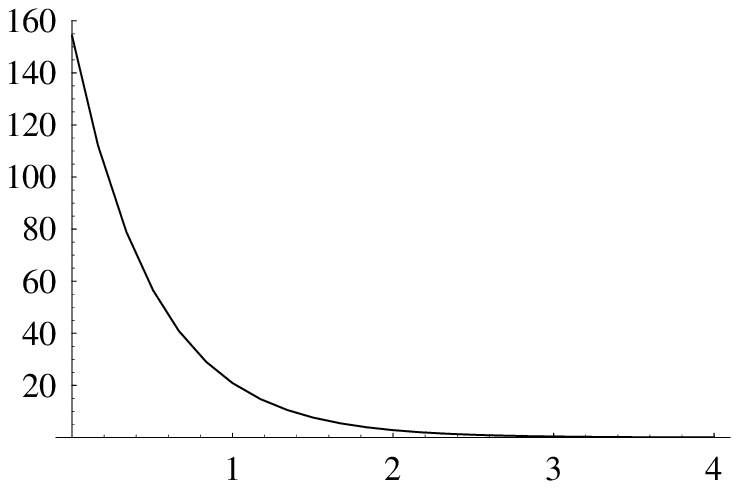}%
\par
{\tiny{
{\textbf{Box 7b.~}
\hbox{$N=5$, $R=0.12$: the function}
\hbox{$\Dd_3(t) := \| u^5(t) \|_3$.
One has $\Dd_3(0) = 154.3...$,} $\Dd_3(0.5) = 57.10...$,
$\Dd_3(1) = 20.88...$, $\Dd_3(1.5) = 7.672... $, $\Dd_3(2) = 2.821...$,
$\Dd_3(4) = 0.05168...$, $\Dd_3(8) =1.733... \times 10^{-5}$,
$\Dd_3(10) = 3.175... \times 10^{-7}$~.
}}
\par}
\label{f7b}
}
}
\hskip 0.4cm
\framebox{
\parbox{2in}{
\includegraphics[
height=1.3in,
width=2.0in
]%
{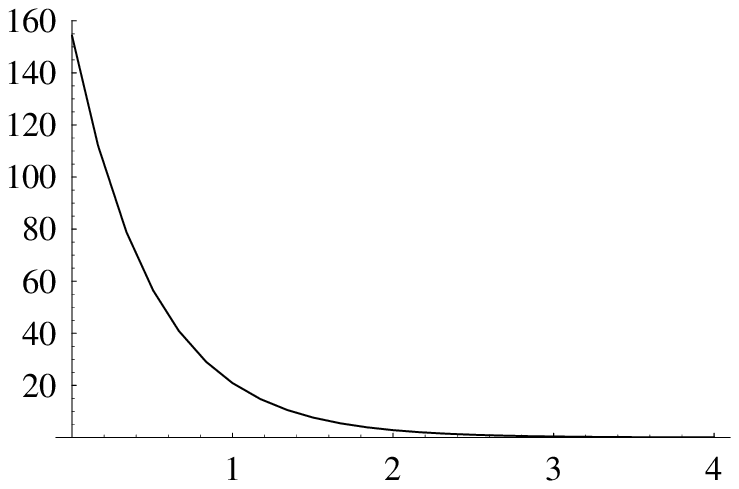}%
\par
{\tiny{
{\textbf{Box 8b.~}
Galerkin method for $R=0.12$: the function
$\Ddd_3(t) := \| \ug(t) \|_3$.
One has $\Ddd_3(0) = 154.3...$, $\Ddd_3(0.5) = 57.10...$,
$\Ddd_3(1) = 20.88...$, $\Ddd_3(1.5) = 7.672... $, $\Ddd_3(2) = 2.821...$,
$\Ddd_3(4) = 0.05168...$, $\Ddd_3(8) =1.736... \times 10^{-5}$,
$\Ddd_3(10) = 3.811... \times 10^{-7}$~.
}}
\par}
\label{f8b}
}
}
\framebox{
\parbox{2in}{
\includegraphics[
height=1.3in,
width=2.0in
]%
{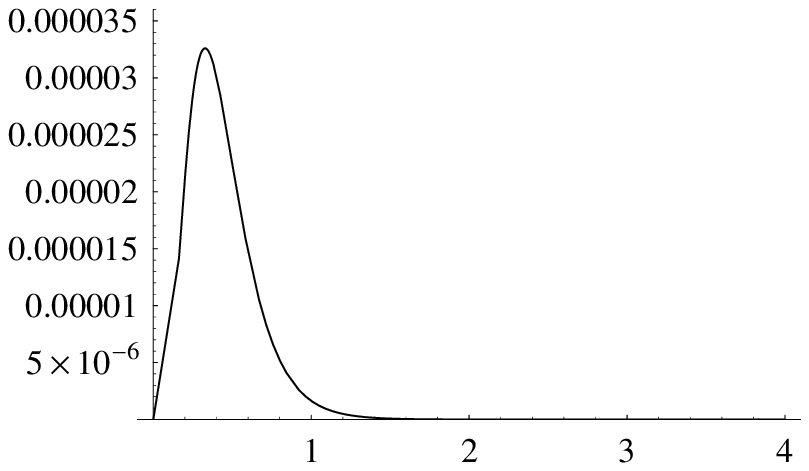}%
\par
{\tiny{
{\textbf{Box 7c.~}
\hbox{$N=5$, $R=0.12$: the function} $\ep_3(t)$.
One has $\ep_3(0) = 0$,
$\ep_3(0.33) = 3.259... \times 10^{-5}$, $\ep_3(1) =1.634...\times 10^{-6} $,
$\ep_3(2) = 4.126... \times 10^{-9}$,
$\ep_3(4) = 1.266... \times 10^{-13}$\, .
}}
\par}
\label{f7c}
}
}
\hskip 0.4cm
\framebox{
\parbox{2in}{
\includegraphics[
height=1.3in,
width=2.0in
]%
{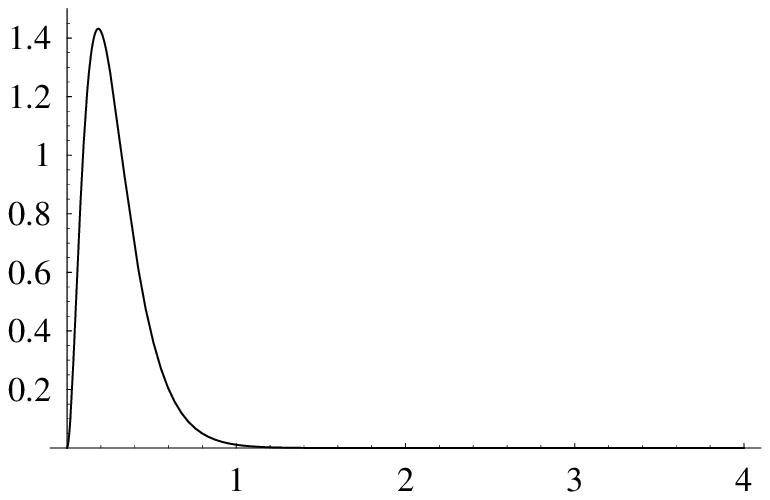}%
\par
{\tiny{
{\textbf{Box 8c.~}
Galerkin method for $R=0.12$:
the function $\eep_3(t)$.
One has $\eep_3(0) = 0$,
$\eep_3(0.18) = 1.431...$,
$\eep_3(1) = 0.01121...$,
$\eep_3(2) = 4.492... \times 10^{-6}$,
$\eep_3(4) = 5.590... \times 10^{-9}$\, .}}
\par}
\label{f8c}
}
}
\framebox{
\parbox{2in}{
\includegraphics[
height=1.3in,
width=2.0in
]%
{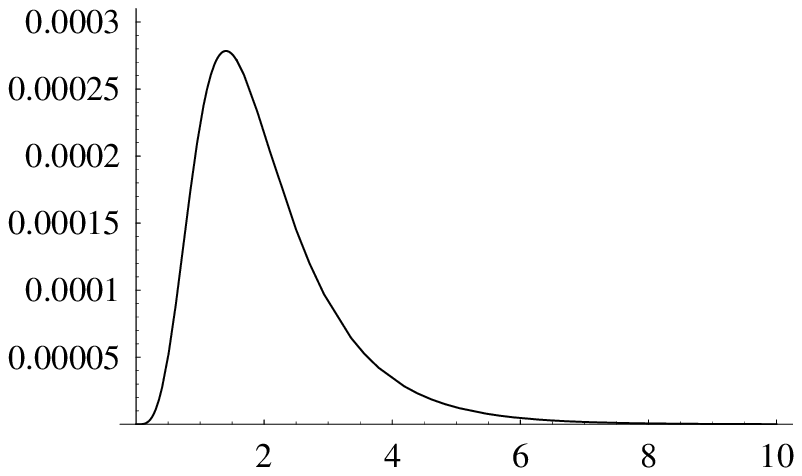}%
\par
{\tiny{
{\textbf{Box 7d.~}
$N=5$, $R=0.12$: the function
$\Rr_3(t)$. This
appears to be \hbox{globally} defined, and vanishing at $+\infty$.
One has $\Rr_3(0) = 0$, $\Rr_3(1) = 2.237... \times 10^{-4}$,
$\Rr_3(1.4) = 2.784... \times 10^{-4}$, $\Rr_3(2) = 2.171... \times 10^{-4}$, $\Rr_3(4) = 3.410... \times 10^{-5}$,
$\Rr_3(10) = 8.431... \times 10^{-8}$~.
}}
\par}
\label{f7d}
}
}
\hskip 3.7cm
\framebox{
\parbox{2in}{
\includegraphics[
height=1.3in,
width=2.0in
]%
{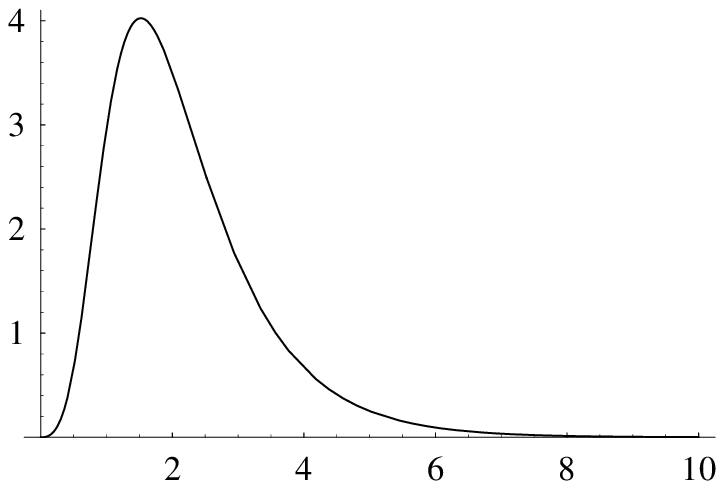}%
\par
{\tiny{
{\textbf{Box 8d.~}
Galerkin method for $R=0.12$:
the function $\Rrr_3(t)$. This appears to
be \hbox{globally} defined, and vanishing at $+\infty$.
One has $\Rrr_3(0)=0$, $\Rrr_3(1)=2.964...$,
$\Rrr_3(1.5)$ $= 4.022...$,
$\Rrr_3(3)$ $= 1.679...$, $\Rrr_3(4)$ $= 0.6665...$,
$\Rrr_3(6)$ $ = 0.09327...$, $\Rrr_3(8)$ $ = 0.01267$,
$\Rrr_3(10) $ $= 1.716 ... \times 10^{-3}$.
}}
\par}
\label{f8d}
}
}
\end{figure}
\vfill \eject \noindent
\vskip 0.7cm \noindent
\textbf{Acknowledgments.} 
This work was partly supported by INdAM, INFN and by MIUR, PRIN 2010
Research Project  ``Geometric and analytic theory of Hamiltonian systems in finite and infinite dimensions''.

\end{document}